\documentclass[a4paper]{article}
\usepackage{lmodern}
\usepackage[T1]{fontenc}
\usepackage{hyperref}
\usepackage{amsfonts,amsmath,amssymb,amsopn}
\usepackage[all]{xy}
\usepackage{vmargin}
\usepackage{makeidx}
\usepackage{graphicx}
\usepackage{tikz}
\usepackage{comment}
\usepackage{titlesec}

\usepackage{amsthm}

\usepackage{color}
\title{Formal confluence of quantum differential operators}
\author{Bernard Le Stum \& Adolfo Quir\'os\thanks{Supported by grant MTM2012-35849 from Ministerio de Econom\'{\i}a y Competitividad (Spain).}}
\date{Version of \today}
\newtheorem{thm}{Theorem}[section]
\newtheorem{prop}[thm] {Proposition}
\newtheorem{cor}[thm] {Corollary}
\newtheorem{lem}[thm] {Lemma}
\theoremstyle{definition}
\newtheorem{dfn}[thm] {Definition}

\newenvironment{xmp}[1][Example]{\begin{trivlist} \item[\hskip \labelsep {\bfseries #1}]}{\end{trivlist}}
 \newenvironment{xmps}[1][Examples]{\begin{trivlist} \item[\hskip \labelsep {\bfseries #1}]}{\end{trivlist}}
\newenvironment{pf}[1][Proof]{\begin{trivlist} \item[\hskip \labelsep {\bfseries #1}]}{\end{trivlist}}
\newenvironment{rmk}[1][Remark]{\begin{trivlist} \item[\hskip \labelsep {\bfseries #1}]}{\end{trivlist}}
 \newenvironment{rmks}[1][Remarks]{\begin{trivlist} \item[\hskip \labelsep {\bfseries #1}]}{\end{trivlist}}
\titleformat{\subsection}[runin]{\normalfont\normalsize\bfseries}{\thesubsection}{1em}{}
\numberwithin{equation}{section}
\begin{document}


\maketitle

\bigskip

\begin{center}
\textbf{Abstract}
\end{center}

\medskip
We prove that a usual differential operator is formally the limit of quantum differential operators.
For this purpose, we introduce the notion of twisted differential operator of infinite level and prove that, formally, such an object is independent of the choice of the twist.
Our method provides explicit formulas.
\medskip

\tableofcontents

\section*{Introduction}

The confluence process consists in replacing differential by finite difference (or $q$-difference) in a differential equation and try to derive solutions of the original equation from solutions of the finite difference (or $q$-difference) equation (see  \cite{LevyLessman92} for example in the finite difference case and \cite{Sauloy00} in the $q$-difference case).
Actually, differential equations and finite difference equations (or $q$-difference equations) may be seen as particular cases of twisted differential equations as explained in \cite{LeStumQuiros15a}. In this more general setting, we are faced with comparing twisted differential equations for different twists.
We will introduce various kinds of twisted differential rings and compare them.
Note that our methods will apply as well to the positive characteristic situation and also to quantum differential equations (with respect to a root of unity).

Let us be more precise. If $R$ is a commutative ring, a twisted $R$-algebra is a commutative $R$-algebra $A$ endowed with an $R$-linear ring endomorphism $\sigma$.
Actually, in \cite{LeStumQuiros15a} we allowed more generally families of endomorphisms satisfying some relations but we will concentrate here on the one dimensional case and we will use either one morphism or a system of roots of one morphism.
One may then consider the $\sigma$-derivations $D$ of $A$, or more generally, the $\sigma$-derivations $D$ of an $A$-module $M$: $R$-linear maps that satisfy the twisted Leibnitz rule
\begin{displaymath}
D(xs) = D(x)s + \sigma(x)D(s).
\end{displaymath}
A $\sigma$-differential module is an $A$-module $M$ endowed with an action by $\sigma$-derivations of the $\sigma$-derivations of $A$.
Assume that $x \in A$ is a $\sigma$-coordinate in the sense for example that there exists a a basis $\partial_{\sigma}$ for the $\sigma$-derivations of $A$ such that $\partial_{\sigma}(x) = 1$ and define the twisted Weyl algebra $\mathrm D_{\sigma}$ as the Ore extension of $A$ by $\sigma$ and $\partial_{\sigma}$ (as in \cite{Bourbaki12}, proposition 1.4).
A $\sigma$-differential module is then the same thing as a $\mathrm D_{\sigma}$-module.
So now, what we want to understand is the relation between $\mathrm D_{\sigma}$ and $\mathrm D_{\tau}$ when $\tau$ is another $R$-linear endomorphism of $A$ with the same twisted coordinate $x$.
In the end, we will be mostly essentially interested in $\tau (x) = x$ (usual case) and $\sigma(x) = qx+h$ (quantum case).

One may usually see a ring of differential operators as some dual of a ring of formal functions.
Doing this directly for $\mathrm D_{\sigma}$ would require to understand the notion of $\sigma$-divided powers on the function side.
This is a very interesting question that we postpone to a forthcoming article.
Here, we will actually replace the twisted Weyl Algebra $\mathrm D_{\sigma}$ with a Grothendieck ring of differential operators $\mathrm D^{(\infty)}_{\sigma}$ (so that the $\sigma$-divided powers live naturally on the differential operator side).
It so happens that the classical construction works incredibly well for this particular generalization.
It is actually sufficient to replace the usual powers of an ideal with the twisted powers introduced in \cite{LeStumQuiros15}.
But now we are faced with two questions: the comparison between $\mathrm D^{(\infty)}_{\sigma}$ and $\mathrm D^{(\infty)}_{\tau}$ on one hand and the comparison of $\mathrm D_{\sigma}$ with $\mathrm D^{(\infty)}_{\sigma}$ on the other.

We will show that, if we denote by $\widehat{\mathrm D}^{(\infty)}_{\sigma}$ the completion along the divided powers of $\partial_{\sigma}$, then there always exists a natural isomorphism $\widehat{\mathrm D}^{(\infty)}_{\sigma} \simeq \widehat{\mathrm D}^{(\infty)}_{\tau}$ (the formal deformation of proposition \ref{dfthm}) and we will be able to give very explicit formulas.
For example, in the case $\tau (x) = x$ and $\sigma(x) = qx$, we can write (over a field of characteristic zero when $q$ is not a root of unity):
\begin{displaymath}
\partial_{q} = \sum_{k=1}^{\infty} \frac 1{k!}(q-1)^{k-1}x^{k-1} \partial^{k} \quad \mathrm{and} \quad \partial = \sum_{k=1}^{\infty} \frac {(1-q)^{k}}{1 -q^k} x^{k-1}\partial_{q}^{k}.
\end{displaymath}

In the quantum situation, when $\sigma(x) = qx +h$ with $q,h \in R$, then there exists a canonical map $\mathrm D_{\sigma} \to \mathrm D^{(\infty)}_{\sigma}$ whose image is exactly the ring of small (or naive) differential operators.
When the $q$-characteristic of $A$, that was introduced in \cite{LeStumQuiros15}, is zero, then all three rings are actually equal and we easily obtain our first confluence theorem (theorem \ref{confalg}) when $R$ is a $\mathbb Q$-algebra: there exists a map $\mathrm D_{\sigma} \to \widehat{\mathrm D}$ with dense image.

But we are actually mainly interested in the positive $q$-characteristic case, which means that $q$ is a primitive $p$-th root of unity.
It is then necessary to use a complete family $\underline q$ of $p^n$-th roots of $q$.
Then we can define what we call a rooted Weyl algebra $\mathrm D_{\underline \sigma}$ by taking the limit on all $\mathrm D_{\sigma_{n}}$ and build a map $\mathrm D_{\underline \sigma} \to \widehat{\mathrm D}$ with dense image.
This is the second confluence theorem (theorem \ref{mainth}).
It is interesting to notice that this theorem puts together the ring $D$ who supports the conjecture of Jacques Dixmier \cite{Dixmier68}  (whose higher version is equivalent to the Jacobian conjecture \cite{BelovKanelKonsevich07}) and an Azumaya algebra where the Morita equivalence could be applied.

We may also notice that  a $p$-adic version of one of our results, the derivation of the formal deformation theorem (proposition \ref{dfthm}) from the density lemma (lemma \ref{density}), would recover Andrea Pulita's realization of a formal solution of a differential module as formal solution of some $\sigma$-module (results in this direction had already obtained by Yves Andr\'e and Lucia Di Vizio in \cite{AndreDiVizio04}). We expect that the quantum refinement could clarify many things and we intend to work this out in the future.

Both authors thank Michel Gros for all the fruitful conversations that we had all together.

For us, a \emph{ring} has an associative multiplication (not commutative in general) and a two-sided unit.
Morphisms of rings are always assumed to preserve the unit.
A module always means a \emph{left} module.
We will essentially consider $1$-twisted rings from \cite{LeStumQuiros15a} and simply call them twisted rings.
Throughout the paper, $R$ will denote a \emph{commutative ring} and $A$ a \emph{twisted commutative $R$-algebra:} a commutative $R$-algebra endowed with an $R$-linear ring endomorphism $\sigma_{A}$.

\section{Twisted principal parts}\label{tpp}

We will introduce here the notion of twisted principal part (functions on twisted infinitesimal neighborhoods of the diagonal).

We will begin by ignoring the endomorphism $\sigma_{A}$ and concentrating on the commutative $R$-algebra $A$.
The tensor product $P_{A/R} := A \otimes_{R} A$ has two $A$-algebra structures, one coming from the action on the left and the other one coming from the action on the right.
Unless otherwise specified, we will use the \emph{left} structure when we see $P_{A/R}$ as an $A$-module.
However, when $M$ is an $A$-module, the notation $P \otimes_{A} M$ will implicitly mean that we use the action of $A$ on the \emph{right} to build the tensor product and that the resulting object will be seen as an $A$-module using the action of $A$ on the \emph{left}.

In practice, we will write $x := x \otimes 1 \in P_{A/R}$ and $\tilde x := 1 \otimes x \in P_{A/R}$.
In other words, with these notations, the action on the left is multiplication by $x$ and the action on the right is multiplication by $\tilde x$.
Any element of $P_{A/R}$ can be written as a finite sum $\sum x_{i} \tilde y_{i}$.
At some point, we will call \emph{Taylor map} the embedding on the right and denote it by
\begin{displaymath}
\xymatrix@R0cm{\theta_{A/R} :& A \ar[r] & P_{A/R}.
\\ & x \ar@{|->}[r] & \tilde x }
\end{displaymath}
We will then call $\tilde x = \theta_{A/R}(x)$ the \emph{Taylor expansion} of $x \in A$ (more on this vocabulary later on).

We will also consider the canonical map corresponding to the projection that forgets the middle term:
\begin{displaymath}
\xymatrix@R0cm{
A \otimes_{R} A \ar[r]^-{\delta_{A/R}} & A \otimes_{R} A \otimes_{R} A
\\
x \otimes y \ar@{|->}[r] & x \otimes 1 \otimes y
}.
\end{displaymath}
It is a morphism of $R$-algebras that may also be seen as a map
\begin{displaymath}
\xymatrix@R0cm{
P_{A/R} \ar[r]^-{\delta_{A/R}} &P_{A/R} \otimes_{A} P_{A/R}
\\
x \tilde y \ar@{|->}[r] & x \otimes \tilde y
},
\end{displaymath}
where $A$ acts on the \emph{right} on the first factor and on the \emph{left} on the second one in the tensor product.

Let $I_{A/R}$ be the kernel of the multiplication morphism
\begin{displaymath}
\xymatrix@R0cm{P_{A/R} \ar[r] & A
\\ x \tilde y \ar@{|->}[r] & xy}
\end{displaymath}
(that corresponds to the diagonal embedding).
The ideal $I_{A/R}$ is generated by the image of the linear map
\begin{displaymath}
\xymatrix@R0cm{A  \ar[r]^-{d_{A/R}} & P_{A/R}
\\ x \ar@{|->}[r] & \xi := \tilde x - x}
\end{displaymath}
(that corresponds to the difference between the projections).
In practice, we will usually drop the index $A/R$ and simply write $P$, $I$, $\theta$, $\delta$ and $d$.

Since we will make regular use of linearization, we introduce it formally now:

\begin{dfn}
Let $M, N$ be two $A$-modules and $u : M \to N$ an $R$-linear map.
Then, the \emph{$A$-linearization} of $u$ is the $A$-linear map
\begin{displaymath}
\xymatrix@R0cm{ P \otimes_{A} M \ar@/^.5cm/[rr] ^-{\tilde u} \ar[r]^{\simeq} &A \otimes_{R} M \ar[r] & N
\\ & x \otimes s \ar@{|->}[r] & xu(s).}
\end{displaymath}
\end{dfn}

Recall that, in the tensor product, we use the action of the \emph{right} on $P$. Therefore, we have
\begin{displaymath}
\forall x,y \in A, \forall s \in M, \quad \tilde u (x\tilde y \otimes s) = xu(ys).
\end{displaymath}

As an example, we see that the multiplication map $P \to A$ is the $A$-linearization of the identity of $A$, seen as an $R$-linear map.

For further use, we also mention the following result, that follows immediately from the definitions:

\begin{lem} \label{compolin}
If we are given two $R$-linear maps $\varphi : M \to N$ and $\psi : L \to M$, then the linearization of $\varphi \circ \psi$ factors as
\begin{equation*} \label{facdif}
\xymatrix@R0cm{ P \otimes_{A} L \ar[r]^-{\delta \otimes \mathrm{Id}} & P \otimes_{A} P \otimes_A L \ar[r]^-{\mathrm{Id} \otimes_{A} \tilde \psi} & P \otimes_{A} M \ar[r]^-{\tilde \varphi} & N.} \quad\Box
\end{equation*}
\end{lem}

Now we make the endomorphism $\sigma_{A}$ enter the game.
We will consider $P := A \otimes_{R} A$ as a twisted $R$-algebra by using $\sigma_{A}$ on the left and the identity the right: in other words, we set $
\sigma_{P} := \sigma_{A} \otimes_{R}  \mathrm{Id}_{A}$.
Alternatively, this is the unique structure of twisted $R$-algebra on $P$ such that 
\begin{displaymath}
\forall x \in A,\quad  \sigma_{P}(x) = \sigma_{A}(x) \quad \mathrm{and} \quad \sigma_{P}(\tilde x) = \tilde x.
\end{displaymath}
In particular, we may also consider $P$ as a twisted $A$-algebra, in the sense that $P$ is endowed with a $\sigma_{A}$-linear ring endomorphism $\sigma_{P}$.
We will often drop the indexes $A$ and $P$ and simply write $\sigma$ for both maps (so that $\sigma(\tilde x) = \tilde x$) when there is no ambiguity.

Before we do anything else, let us prove the following result, which is quite elementary, but very useful:

\begin{lem} \label{transl}
If $x \in A$ and $\xi :=\tilde x - x \in P_{A/R}$, then,
\begin{displaymath}
\sigma_{P}(\xi) = \xi + y \quad \mathrm{with} \quad y:= x - \sigma_{A}(x).
\end{displaymath}
\end{lem}

\begin{pf}
We have
\begin{displaymath}
\sigma(\xi) = \sigma(\tilde x - x) = \sigma(\tilde x) - \sigma(x) = \tilde x - \sigma(x) = \tilde x - x + x - \sigma(x) = \xi + y. \quad \Box
\end{displaymath}
\end{pf}

\begin{lem} \label{kersig} The kernel of the $A$-linearization
\begin{displaymath}
\xymatrix@R0cm{ P_{A/R} \ar[r]^-{\tilde \sigma_{A}} & A
\\ x \tilde y \ar@{|->}[r] & x\sigma_{A}(y)}
\end{displaymath}
of $\sigma_{A}$ is $\sigma_{P}(I_{A/R})$.
\end{lem}

Recall that $\sigma(I)$ denotes the \emph{ideal} generated by the image of $I$.
As a consequence, since $I$ is generated by the image of $d$, we see that $\ker \tilde \sigma$ will be generated by the image of
\begin{displaymath}
\xymatrix@R0cm{A  \ar[r]^-{\sigma \circ d} & P
\\ x \ar@{|->}[r] & \tilde x - \sigma(x)}.
\end{displaymath}

\begin{pf}
First of all, we have
\begin{displaymath}
\forall x \in A, \quad \tilde \sigma( \tilde x - \sigma(x)) = \sigma(x) - \sigma(x) = 0,
\end{displaymath}
and it follows that $\sigma(I) \subset \ker \tilde \sigma$.
Conversely, by definition, we have
\begin{displaymath}
\forall x \in A, \quad \sigma(x)  \equiv \tilde x \mod \sigma(I).
\end{displaymath}
Therefore, if $f := \sum x_{i} \tilde y_{i} \in \ker \tilde \sigma$, we have
\begin{displaymath}
f \equiv  \sum x_{i}\sigma(y_{i}) = \tilde \sigma (f) = 0 \mod \sigma(I),
\end{displaymath}
and we see that $\ker \tilde \sigma \subset \sigma(I)$.
$\quad \Box$
\end{pf}

\begin{rmks}
\begin{enumerate}
\item
It is sometimes convenient to use the bimodule language.
An \emph{$A$-sesquimodule} $M$ is an $A$-bimodule such that
\begin{displaymath}
\forall x \in A, s \in M, \quad \sigma_{A}(x) \cdot s = s \cdot x.
\end{displaymath}
Note that we are using the \emph{reverse} convention from Andr\'e's in \cite{Andre01} so that forgetting the right action induces an equivalence (an isomorphism) between $A$-sesquimodules and left $A$-modules (we will use this identification).
\item 
The $R$-algebra $P$ has a canonical $A$-bimodule structure which is completely independent of the choice of $\sigma_{A}$.
If we endow $A$ with its sesquimodule structure, then the linearization $\tilde \sigma$ of $\sigma_{A}$ is a morphism of $A$-bimodules: we always have 
\begin{displaymath}
\forall x, y \in A, f \in P, \quad \tilde\sigma(x\cdot f\cdot y)=\tilde \sigma(xf\tilde y) = x\tilde \sigma(f)\sigma(y)=x\cdot \tilde\sigma(f) \cdot y.
\end{displaymath}
It follows that $\ker \tilde \sigma$ has a natural structure of $A$-bimodule.
Actually, since $\tilde \sigma$ is a ring homomorphism, then $\ker \tilde \sigma$ is an ideal and therefore automatically an $A$-bimodule.
\end{enumerate}
\end{rmks}

Recall from \cite{LeStumQuiros15} that the \emph{twisted powers}  of $I$ are
\begin{displaymath}
I^{(0)} = P, \quad I^{(1)} = I , \quad I^{(2)} := I\sigma(I), \quad \ldots, \quad  I^{(n)} := I\sigma(I) \cdots \sigma^{n-1}(I),
\end{displaymath}
where images and products of ideals are meant as \emph{ideals}.
We will write $I^{(n)_{\sigma}}_{A/R}$ when we want to make clear the dependence on $\sigma$ and $A/R$.

\begin{dfn}
The $A$-module of \emph{twisted principal parts of order $n \in \mathbb N$ (and infinite level)} of $A$ is
\begin{displaymath}
P_{A/R,(n)_{\sigma}} : = P_{A/R}/I_{A/R}^{(n+1)_{\sigma}}.
\end{displaymath}
The $A$-module of \emph{twisted principal parts of infinite order (and infinite level)} of $A$ is the \emph{twisted completion}:
\begin{displaymath}
\widehat P_{A/R,\sigma} := \varprojlim P_{A/R,(n)_{\sigma}}.
\end{displaymath}
\end{dfn}

Note that these $A$-modules all have a natural structure of $R$-algebra and that the definition also makes sense for $n = -1$ so that $P_{A/R,(-1)_{\sigma}} = 0$.
Again, we will often drop the indexes $A/R$ when we believe that there is no risk of confusion and simply write $P_{(n)_{\sigma}}$ and $\widehat P_{\sigma}$.
We will also drop the index $\sigma$ when $\sigma_{A} = \mathrm{Id}_{A}$.

\begin{xmps}
\begin{enumerate}
\item When $A$ is trivially twisted, which means that $\sigma_{A} := \mathrm{Id}_{A}$, this notion of principal parts coincides with the usual one (definition 16.3.1 of \cite{EGA44}), and therefore many of the basic objects we will construct are twisted versions of those in  \cite{EGA44}.
\item When $A = R[x]$, we have $P = R[x, \xi]$ with $\xi = \tilde x - x$ and $I = (\xi)$.
Moreover, $\sigma(\xi) = \xi + y$ with $y = x - \sigma(x)$.
It follows that
\begin{displaymath}
\sigma^n(\xi) = \xi + (n)_{\sigma}(y) = \xi + x - \sigma^{n}(x),
\end{displaymath}
with $(n)_{\sigma} := 1 + \cdots + \sigma^{n-1}$.
Therefore, we have
\begin{displaymath}
P_{(n)} = R[x, \xi]/\prod_{i=0}^n \left(\xi + (i)_{\sigma}(y)\right) =  R[x, \xi]/\prod_{i=0}^n \left(\xi + x - \sigma^{n}(x)\right).
\end{displaymath}
\begin{enumerate}
\item In the case $\sigma(x) = x $, we get $P_{(n)} = R[x, \xi]/\xi^{n+1}$ as expected.
\item More generally, if we assume that $\sigma(x) = x + h$ with $h \in R$, we obtain
\begin{displaymath}
P_{(n)} = R[x, \xi]/\prod_{i=0}^n \left(\xi -ih\right).
\end{displaymath}
\item On the other hand, if we let $\sigma(x) = qx$ with $q \in R$, we find
\begin{displaymath}
P_{(n)} = R[x, \xi]/\prod_{i=0}^n \left( \xi + (1-q^i)x\right).
\end{displaymath}
\end{enumerate}
\end{enumerate}
\end{xmps}

When $R \to R'$ is a homomorphism of commutative rings, we  endow $A' := R' \otimes_R A$ with $\sigma_{A'} := \mathrm{Id}_{R'} \otimes_{R} \sigma_{A}$.

\begin{prop} \label{extn} 
Let $R \to R'$ be a homomorphism of commutative rings and $A' := R' \otimes_R A$.
Then, we have for all $n \in \mathbb N$,
\begin{displaymath}
A' \otimes_A P_{A/R,(n)_{\sigma}} \simeq P_{A'/R',(n)_{\sigma}}.
\end{displaymath}
\end{prop}

\begin{pf}
If we let $P' := P_{A'/R'}$, then there exists an canonical isomorphism $A' \otimes_{A} P \simeq R' \otimes_{R} P \simeq P'$.
Moreover, if we denote by $I'$ the kernel of multiplication on $P'$, we have $A' \otimes_{A} I \simeq I'$.
And finally, $\mathrm{Id}_{A} \otimes_{A} \sigma_{P}$ corresponds to $\sigma_{P'}$ under this isomorphism.
Our assertion is therefore an immediate consequence of right exactness of tensor product.
$\quad \Box$
\end{pf}

Recall from \cite{LeStumQuiros15} that a twisted $A$-algebra is an $A$-algebra $B$ endowed with a $\sigma_{A}$-linear ring endomorphism $\sigma_{B}$.

\begin{prop} \label{locord}
If $B$ is a twisted commutative $A$-algebra, then there exists a canonical morphism of $B$-algebras
\begin{displaymath}
B \otimes_A P_{A/R,(n)_{\sigma}} \to P_{B/R,(n)_{\sigma}}.
\end{displaymath}
When $B$ is a quotient (resp. a localization) of $A$ this map is surjective (resp. bijective).
\end{prop}
 
Recall from Definition 1.7 of \cite{LeStumQuiros15a} that we call such a $B$ a \emph{twisted quotient} (resp. \emph{twisted localization}) of $A$.

\begin{pf}
The morphism of twisted $R$-algebras $A \to B$ extends naturally to a morphism of twisted $R$-algebras $P_{A} \to P_{B}$.
Since $I_{A}$ is sent into $I_{B}$, we see that, for all $n \in \mathbb N$, $\sigma^n(I_{A})$ is sent into $\sigma^n(I_{B})$ and the first assertion formally follows.
In the case of a quotient map, all the maps involved are surjective.

Now, if $B := S^{-1}A$ is a twisted localization of $A$, then $P_{B}$ is the localization of $P_{A}$ with respect to the monoid $S'$ generated by $S$ and $\tilde S$, and we have $I_{B} = P_{B} \otimes_{P_{A}} I_{A}$.
It immediately follows that for all $n \in \mathbb N$, we have $\sigma^n(I_{B}) = P_{B} \otimes_{P_{A}} \sigma^n(I_{A})$, and therefore $I_{B}^{(n)} = P_{B} \otimes_{P_{A}} I_{A}^{(n)}$.
Thus we see that
\begin{displaymath}
P_{B,(n)_{\sigma}} = P_{B}/(P_{B} \otimes_{P_{A}} I_{A}^{(n+1)}) =  P_{B} \otimes_{P_{A}} P_{A,(n)} = B \otimes_A P_{A,(n)} \otimes_{A} B.
\end{displaymath}
We need to remove the $B$ on the right hand side and it is sufficient to show that $\tilde x$ is invertible in $B \otimes_A P_{A,(n)}$ whenever $x \in S$.
But we have
\begin{displaymath}
\prod_{i=0}^n (\tilde x - \sigma^i(x)) = \prod_{i=0}^n \sigma^i(\tilde x - x) \in I_{A}^{(n+1)},
\end{displaymath}
from which we derive that there exists $f \in P_{A}$ such that
\begin{displaymath}
f\tilde x \equiv \prod_{i=0}^n \sigma^i(x) \mod I_{A}^{(n+1)}.
\end{displaymath}
Since $B = S^{-1} A$, we must have $\sigma(S) \subset B^\times$ and it follows that $\prod_{i=0}^n \sigma^i(x) \in B^\times$.
Thus, we see that $f\tilde x$ is invertible in $B \otimes_A P_{A,(n)}$ and it follows that $\tilde x$ is invertible too.
$\quad \Box$
\end{pf}

As an illustration, we can give explicit formulas in the quantum situation.
Recall that we introduced in \cite{LeStumQuiros15} the notion of \emph{twisted powers} of an element in a twisted ring.
In particular, for $f \in P$, we will have
\begin{displaymath}
f^{(0)} := 1, \quad f^{(1)} := f, \quad f^{(2)} = f\sigma(f), \quad \ldots, \quad f^{(n+1)} = f\sigma(f) \cdots \sigma^n(f).
\end{displaymath}
Recall also that the \emph{quantum binomial coefficients} are defined by induction (see \cite{LeStumQuiros15} for example) as
\begin{displaymath}
{n \choose k}_{q} := {n-1 \choose k-1}_{q} + q^k{n-1 \choose k}_{q}.
\end{displaymath}

\begin{prop}
Assume $\sigma(x) = qx$ with $q \in R$ and let $\xi = \tilde x -x$.
Then, we have
\begin{equation} \label{binxi}
\forall n \in \mathbb N, \quad \xi^{(n)} = \sum_{j=0}^n (-1)^{j}  {n \choose j}_q q^{\frac {j(j-1)}2}x^{j}\tilde x^{n-j}
\end{equation}
and 
\begin{displaymath}
\forall n \in \mathbb N, \quad \tilde x^n = \sum_{i=0}^n {n \choose i}_{q} x^i\xi^{(n-i)}.
\end{displaymath}
\end{prop}

\begin{pf}
The first equality is essentially the quantum binomial formula (see proposition 2.14 in \cite{LeStumQuiros15}):
\begin{displaymath}
(\tilde x - x)^{(n)} = \sum_{j=0}^n  {n \choose j}_q (-x)^{(j)}\tilde x^{(n-j)}.
\end{displaymath}

For the second one, we compute the right hand side with the help of formula \eqref{binxi}:
\begin{eqnarray*}
\sum_{i=0}^n {n \choose i}_{q} x^i\xi^{(n-i)}
&=& \sum_{i=0}^n {n \choose i}_{q} x^i \left(\sum_{j=0}^{n-i} (-1)^j {n-i \choose j}_{q} q^{\frac {j(j-1)}2} x^j\tilde x^{n-i-j}\right)
\\
&=& \sum_{k=0}^{n} \left( \sum_{i=0}^k  {n \choose i}_{q}  {n-i \choose k-i}_{q} (-1)^{k-i}  q^{\frac {(k-i)(k-i-1)}2} \right) x^{k}\tilde x^{n-k}
\end{eqnarray*}
after rewriting $i+j=k$.
Now, we have (using corollaries 2.7 and 2.8 in \cite{LeStumQuiros15}, for example)
\begin{displaymath}
\sum_{i=0}^k  {n \choose i}_{q}  {n-i \choose k-i}_{q} (-1)^{k-i}  q^{\frac {(k-i)(k-i-1)}2} = {n \choose k}_{q} \sum_{i=0}^k  {k \choose i}_{q}  (-1)^{k-i}  q^{\frac {(k-i)(k-i-1)}2}.
\end{displaymath}
But the quantum binomial formula again implies that for $k > 0$, we have
\begin{displaymath}
 \sum_{i=0}^k  {k \choose i}_{q}  (-1)^{k-i}  q^{\frac {(k-i)(k-i-1)}2} = \prod_{i=0}^{k-1} (1-q^i) = 0.
\end{displaymath}
And it follows that
\begin{displaymath}
\sum_{i=0}^n {n \choose i}_{q} x^i\xi^{(n-i)} = \tilde x^n,
\end{displaymath}
as asserted.
$\quad \Box$
\end{pf}

We will need a slightly stronger notion of coordinate than the one used in \cite{LeStumQuiros15a}: 

\begin{dfn} \label{twistdef}
Let $x \in A$ and $\xi :=  \tilde x - x \in P_{A/R}$.
\begin{enumerate}
\item 
Then, $x$ is a \emph{twisted coordinate} for $A$ over $R$ if for all $n \in \mathbb N$, $P_{A/R,(n)_{\sigma}}$ is freely generated as an $A$-module by the images of $1, \xi, \xi^{2}, \ldots, \xi^{n}$.
\item
If $x$ is a twisted coordinate such that $\sigma(x) = qx + h$ with $q, h \in R$, we will call it a \emph{quantum coordinate} or \emph{$q$-coordinate} and call $A$ a \emph{quantum $R$-algebra} or \emph{$q$-$R$-algebra}.
\end{enumerate}
\end{dfn}

\begin{xmps}
\begin{enumerate}
\item
When $A/R$ is smooth (of pure relative dimension one) and $\sigma_{A} = \mathrm{Id}_{A}$, then a twisted coordinate is nothing but an \emph{\'etale coordinate}: it means that the map
\begin{displaymath}
\xymatrix@R0cm{R[T]  \ar[r] & A
\\ T \ar@{|->}[r] & x}
\end{displaymath}
is \'etale.
\item If $A = R[x]$, then $x$ is always a twisted coordinate, whatever $\sigma$ is.
\end{enumerate}
\end{xmps}

\begin{prop} \label{locco}
\begin{enumerate}
\item If $R \to R'$ is a homomorphism of commutative rings and $x$ is a twisted coordinate on $A$, then $x$ becomes a twisted coordinate on $A' := R' \otimes_{R} A$ (relatively to $R'$).
\item If $B$ is a twisted localization of $A$ and $x$ is a twisted coordinate on $A$, then $x$ becomes a twisted coordinate on $B$.
\end{enumerate}
\end{prop}

\begin{pf}
Follows from propositions \ref{extn} and \ref{locord}.
$\quad \Box$
\end{pf}

\begin{rmk}
Assume $A$ is a quantum $R$-algebra so that there exists a twisted coordinate $x$ on $A$ and $q,h \in R$ such that $\sigma(x) = qx+h$.
Let us still denote by the same letter $x$ an indeterminate over $R$ and by $\sigma$ again the endomorphism of $R[x]$ given by the same formula.
Then, $A$ becomes an $R[x]$-twisted algebra and we have a canonical isomorphism (compare basis on both sides):
\begin{displaymath}
A \otimes_{R[x]} P_{R[x]/R,(n)_{\sigma}} \simeq P_{A/R,(n)_{\sigma}}.
\end{displaymath}
\end{rmk}

In the next statement, we use the letter $\xi$ as an indeterminate over $A$ so that $A[\xi]$ denotes the polynomial ring and $A[\xi]_{\leq n}$ the submodule of polynomial of degree at most $n$.
Ultimately, this should not create any confusion due to corollary \ref{insid} below.

\begin{prop} \label{sigis}
Let $x \in A$ and $y := x - \sigma(x)$.
We endow $A[\xi]$ with the unique $\sigma_{A}$-linear endomorphism such that $\sigma(\xi) = \xi + y$.
Then, $x$ is a twisted coordinate on $A$ if and only if the morphism of twisted algebras
\begin{displaymath}
\xymatrix@R0cm{ \phi\,:\, A[\xi] \ar[r]& P_{A/R}
\\ \xi \ar@{|->}[r] & \tilde x - x}
\end{displaymath}
induces for all $n \in \mathbb N$, an isomorphism of $A$-algebras $A[\xi]/\xi^{(n+1)} \simeq P_{A/R,(n)_{\sigma}}$.
\end{prop}

\begin{pf}
First of all, it follows from lemma \ref{transl} that there exists such a morphism for all $n \in \mathbb N$.
On the other hand,
\begin{displaymath}
\xi^{(n+1)}= \prod_0^{n} (\xi + (x - \sigma^i(x)))
\end{displaymath}
is a monic polynomial of degree $n+1$.
Then, euclidean division tells us that the composite map
\begin{displaymath}
\xymatrix@R0cm{ A[\xi]_{\leq n} \ar[r] & A[\xi] \ar[r] & A[\xi]/  \xi^{(n+1)}}
\end{displaymath}
is an isomorphism of $A$-modules.
Therefore the condition on $\phi$ is equivalent to the fact that the map
\begin{displaymath}
\xymatrix@R0cm{  A[\xi]_{\leq n} \ar[r] & P_{(n)_{\sigma}}
\\ \xi \ar@{|->}[r] & \overline{\tilde x - x}}
\end{displaymath}
is bijective.
And this exactly means that $P_{(n)_{\sigma}}$ is freely generated by the $n+1$ first powers of the images of $\tilde x - x$.
$\quad \Box$
\end{pf}

When the polynomial ring $A[\xi]$ is endowed with a structure of $\sigma_{A}$-algebra, we will set
\begin{displaymath}
A[[\xi]]_{\sigma} := \varprojlim A[\xi]/  \xi^{(n+1)}.
\end{displaymath}

\begin{cor} \label{insid}
With the same hypothesis, $x$ is a twisted coordinate on $A$ if and only if there exists an isomorphism of $A$-algebras
\begin{displaymath}
\xymatrix@R0cm{  A[[\xi]]_{\sigma} \ar[r]^-{\simeq} &\widehat P_{A/R,\sigma}\,.
\\ \xi \ar@{|->}[r] & \tilde x - x &  \Box}
\end{displaymath}
\end{cor}

\begin{cor} \label{free}
Let $x \in A$ and $\xi :=  \tilde x - x \in P$.
Then, the following conditions are equivalent:
\begin{enumerate}
\item $x$ is a twisted coordinate on $A$
\item for all $n \in \mathbb N$, the $A$-module $P_{(n)_{\sigma}}$ is freely generated by the images of $1, \xi, \xi^{(2)}, \ldots, \xi^{(n)}$
\item for all $n \in \mathbb N$, the $A$-module $I^{(n)}/I^{(n+1)}$ is free of rank one on the image of $\xi^{(n)}$.
$\quad \Box$
\end{enumerate}
\end{cor}

\begin{cor}
If $A$ is a twisted localization of $R[x]$, then $x$ is a twisted coordinate on $A$.
\end{cor}

\begin{pf}
Using the second part of proposition \ref{locco}, we may assume that $A = R[x]$ in which case, this is a trivial consequence of proposition \ref{sigis}.
$\quad \Box$
\end{pf}

\section{Twisted differential forms}

In this section, we study the module of twisted differential forms (of degree one) and make the link with twisted derivations.
We use same notations as before.

\begin{dfn} The $A$-module of \emph{twisted differential forms} on $A/R$ is
\begin{displaymath}
\Omega_{A/R, \sigma}^1 := I_{A/R}/I_{A/R}^{(2)_{\sigma}}.
\end{displaymath}
\end{dfn}

Again, we will often drop the index $A/R$.
Since we implicitly endow $P$ with the action of $A$ on the left, we will also always see $\Omega^1_{\sigma}$ as an $A$-module through the action on the left.

\begin{xmps}
\begin{enumerate}
\item
When $A$ is trivially twisted, then $\Omega^1_{\sigma} = I/I^2$ is the usual module of differential forms of $A$ over $R$.
\item
If $A = R[x]$ is endowed with any $R$-algebra endomorphism $\sigma$, then $\Omega^1_{\sigma}$ is free of rank $1$:
with the notations of lemma \ref{transl}, we have
\begin{displaymath}
\Omega^1_{\sigma} \simeq \xi R[x,\xi]/\xi(\xi + y) \simeq R[x,\xi]/(\xi + y) \simeq R[x].
\end{displaymath}
\end{enumerate}
\end{xmps}

\begin{rmks}
\begin{enumerate}
\item
Clearly, $\Omega^1_{\sigma}$ has a natural structure of $A$-bimodule as a quotient of two ideals of $P$.
It happens that this is identical to its $A$-sesquimodule structure:
by definition, if $x \in A$, then $\tilde x \equiv \sigma(x) \mod \sigma(I)$ and it follows that
\begin{equation}\label{modisig}
\forall f \in I, \quad f\tilde x \equiv f\sigma(x) \mod I\sigma(I).
\end{equation}
\item The $\Omega^1_{\sigma}$ that appears in proposition 1.4.2.1 of \cite{Andre01} is exactly the same as ours (Andr\'e calls $k$ what we call $R$).
\item Formula \ref{modisig} is exactly the first step of the braiding described by Max Karoubi and Mariano Su\'arez-\'Alvarez  in \cite{KaroubiSuarez03}.
\item One can define more generally the \emph{twisted de Rham complex} $\Omega^\bullet_{\sigma}$ of $A$ as the quotient of the non commutative tensor algebra of $I$ by the graded differential ideal generated by $I\sigma(I)$.
We will not consider this complex here.
\end{enumerate}
\end{rmks}

\begin{prop} \label{splitex}
There exists a split exact sequence
\begin{displaymath}
\xymatrix{0 \ar[r] &  \Omega^1_{A/R,\sigma} \ar[r] & P_{A/R,(1)_{\sigma}} \ar[r]^-{\tilde \sigma_{A}} & A \ar[r] & 0}.
\end{displaymath} 
\end{prop}

\begin{pf}
There exists such an exact sequence by definition of $P_{(1)_{\sigma}}$ and $\Omega^1_{\sigma}$.
The $A$-module structure of $P$ provides a section of $\widetilde \sigma$.
$\quad \Box$
\end{pf}

\begin{prop} \label{omegfonc1} \label{omegfonc2}
\begin{enumerate}
\item  If $R \to R'$ is a  homomorphism of commutative rings and $A' := R' \otimes_R A$, then there exists an isomorphism
\begin{displaymath}
A' \otimes_A \Omega_{A/R, \sigma}^1 \simeq \Omega_{A'/R',  \sigma}^1.
\end{displaymath}
\item
If $B$ is a twisted commutative $A$-algebra, then there exists a canonical $B$-linear map
\begin{equation} \label{compom}
B \otimes_{A} \Omega_{A/R,  \sigma}^1 \to \Omega_{B/R, \sigma}^1.
\end{equation}
When $B$ is a quotient (resp. a localization) of $A$ this map is surjective (resp. bijective).
\end{enumerate}
\end{prop}

\begin{pf} Using proposition \ref{splitex}, this follows from propositions \ref{extn}  and \ref{locord}.
\end{pf}

\begin{rmk}
This last result does \emph{not} hold however if we only require $A \to B$ to be an \'etale map (and not a localization map) as the following example shows.
Let $R$ be any field of characteristic different from 2, $A := R [x]$ with $\sigma_{A} := \mathrm{Id}_{A}$ and $B :=  R [x,x^{-1}]$ with $\sigma_{B}(x) := -x$.
Then the morphism $x \mapsto x^2$ is an \'etale twisted morphism but the morphism \eqref{compom} is the zero map.
More precisely, if $\xi = \tilde x - x$, we have 
\begin{displaymath} \label{compom2}
B \otimes_A \Omega_{A,\sigma}^1 = (\xi)/(\xi^2)  \quad \mathrm{and} \quad \Omega_{B,\sigma}^1 =  (\xi)/(\xi^2 +2x\xi) ,
\end{displaymath}
where the ideals are taken inside $R[x,x^{-1},\xi]$, and morphism \eqref{compom2} is induced by $\xi \mapsto \xi^2 + 2x\xi$.
\end{rmk}

Recall from \cite{LeStumQuiros15a} that a \emph{twisted derivation} of $A$ is an $R$-linear map into an $A$-module $M$ that satisfies the \emph{twisted Leibnitz rule}:
\begin{displaymath}
\forall x,y \in A, \quad D(xy) = x D(y) + \sigma(y)D(x).
\end{displaymath}
They form an $A$-module $\mathrm {Der}_{R,\sigma}(A,M)$.

\begin{prop} \label{deriv} The canonical map $A \to \Omega^1_{A/R,\sigma}$ induced by $d$ is a twisted derivation.
It provides us with a natural isomorphism
\begin{equation}\label{stan}
\xymatrix@R0cm{ \mathrm {Hom}_{A}(\Omega^1_{A/R,\sigma}, M)  \ar[r]^-\simeq & \mathrm {Der}_{R,\sigma}(A,M)
\\ u  \ar@{|->}[r] & D := u \circ d}
\end{equation}
whenever $M$ is an $A$-module.
\end{prop}

In the future, we will also denote by $d : A \to \Omega^1_{\sigma}$ this universal twisted derivation when there is no risk of confusion.

\begin{pf}
Using formula \eqref{modisig}, we see that, inside $P$, we have
\begin{eqnarray*}
yd(x) + \sigma(x)d(y) 
&=& y(\tilde x - x) + \sigma(x)(\tilde y -y)
\\
&\equiv& y(\tilde x - x) + \tilde x(\tilde y -y) = \widetilde {xy} - xy = d(xy) \mod I\sigma(I).
\end{eqnarray*}
It follows that the induced map $d : A \to \Omega^1_{\sigma}$ is indeed a twisted derivation.
This also implies that the map in \eqref{stan} is well defined.
And it is clearly injective because $I$ is generated by the image of $d$.

We now show that it is surjective.
If $D$ is a twisted derivation of $M$, we can consider its linearization
\begin{displaymath}
\xymatrix@R0cm{ P \ar[r]^{\tilde D} & M
\\  x \tilde y \ar@{|->}[r] & xD(y)}.
\end{displaymath}
By definition, the ideal $I\sigma(I)$ is generated by elements of the form
\begin{displaymath}
f = (\tilde x - x)(\tilde y - \sigma(y)) = \tilde x \tilde y - x \tilde y - \sigma(y) \tilde x + x\sigma(y),
\end{displaymath}
and we have
\begin{displaymath}
\tilde D(f) = D(xy) - xD(y) - \sigma(y)D(x) + x\sigma(y) D(1) = 0
\end{displaymath}
because $D$ is a twisted derivation (and in particular $D(1) = 0$).
It follows that $\tilde D$ factors through $P/I\sigma(I)$ and we may consider the induced map $u : \Omega^1_{\sigma} \to M$.
It only remains to notice that we have
\begin{displaymath}
\forall x \in A, \quad u(d(x)) = \tilde D(d(x))
= \tilde D (\tilde x - x) = D (x) - x D(1) = D(x).
\quad \Box
\end{displaymath}
\end{pf}

\begin{rmk}
There exists a very elegant proof of this last result through the theory of bimodules.
It is based on the fact (see proposition 17 in \cite{Bourbaki70}, Chapter III, section 10) that $I$ is universal for bimodule derivations: there exists a natural isomorphism
\begin{displaymath}
\xymatrix@R0cm{ \mathrm {Hom}_{A\mathrm{-Bim}}(I, M)  \ar[r]^-\simeq & \mathrm {Der}_{R}(A,M)
\\ u  \ar@{|->}[r] & D := u \circ d},
\end{displaymath}
where the right-hand side stands for bimodule derivations (see proposition 1.4.2.1 of \cite{Andre01}).
\end{rmk}

As an immediate consequence of the proposition, if we write $T_{A/R,\sigma} := \mathrm {Der}_{R,\sigma}(A,A)$, that we will often abbreviate to $T_{\sigma}$, we obtain the following:

\begin{cor} \label{dualid}  The $A$-module $T_{A/R,\sigma}$ is the dual of $\Omega^1_{A/R,\sigma}$.
$\quad \Box$
\end{cor}

\begin{prop} \label{corfonc} 
Assume that $\Omega^1_{A/R,\sigma}$ is projective of finite rank.
Then, if $M$ is an $A$-module, we have the following:
\begin{enumerate}
\item If $R \to R'$ is a base extension and $A' = R' \otimes_{R} A$, then there exists a canonical isomorphism
\begin{displaymath}
R' \otimes_{R} \mathrm {Der}_{R,\sigma}(A,M) \simeq  \mathrm {Der}_{R',\sigma}(A',A' \otimes_{A} M).
\end{displaymath}
\item If $B$ is a twisted $A$-algebra, there exists a canonical map
\begin{displaymath}
B \otimes_{A} \mathrm {Der}_{R,\sigma}(A,M) \leftarrow \mathrm {Der}_{R,\sigma}(B,B \otimes_{A} M).
\end{displaymath}
It is injective (resp. bijective) when $B$ is a quotient (resp. a localization) of $A$.
\end{enumerate}
\end{prop}

\begin{pf}
Both assertions follow from proposition \ref{omegfonc1}.
More precisely, in the first case, we have
\begin{displaymath}
R' \otimes_{R} \mathrm {Der}_{R,\sigma}(A,M) \simeq R' \otimes_{R} \mathrm {Hom}_{A}(\Omega^1_{A/R,\sigma}, M) \simeq \mathrm {Hom}_{A}(\Omega^1_{A/R,\sigma}, A' \otimes_{A} M) 
\end{displaymath}
because $\Omega^1_{A/R,\sigma}$ is projective of finite rank, and then
\begin{eqnarray*}
\mathrm {Hom}_{A}(\Omega^1_{A/R,\sigma}, A' \otimes_{A} M) 
&\simeq& \mathrm {Hom}_{A'}(A' \otimes_{A} \Omega^1_{A/R,\sigma}, A' \otimes_{A} M)
\\
&\simeq& \mathrm {Hom}_{A'}(\Omega^1_{A'/R',\sigma}, A' \otimes_{A} M) \simeq \mathrm {Der}_{R',\sigma}(A',A' \otimes_{A} M).
\end{eqnarray*}
The proof of the second assertion follows exactly the same lines with the same arguments.
We have
\begin{displaymath}
B \otimes_{A} \mathrm {Der}_{R,\sigma}(A,M) \simeq B \otimes_{A} \mathrm {Hom}_{A}(\Omega^1_{A/R,\sigma}, M) \simeq \mathrm {Hom}_{A}(\Omega^1_{A/R,\sigma}, B \otimes_{A} M) \simeq
\end{displaymath}
\begin{displaymath}
 \mathrm {Hom}_{B}(B \otimes_{A} \Omega^1_{A/R,\sigma}, B \otimes_{A} M)
\leftarrow \mathrm {Hom}_{A}(\Omega^1_{B/R,\sigma}, B \otimes_{A} M) \simeq \mathrm {Der}_{R,\sigma}(B',B \otimes_{A} M).
\end{displaymath}
And the only map which is not always an isomorphism will be injective (resp. bijective) when $B$ is a quotient (resp. a localization) of $A$.
$\quad \Box$
\end{pf}

The following immediate consequence is worth stating:

\begin{cor} 
If $\Omega^1_{A/R,\sigma}$ is projective of finite rank, then $\mathrm T_{A/R,\sigma}$ and $\Omega^1_{A/R,\sigma}$ are dual to each other and we have
\begin{enumerate}
\item If $R \to R'$ is a base extension and $A' = R' \otimes_{R} A$, then there exists a canonical isomorphism
\begin{displaymath}
R' \otimes_{R} \mathrm T_{A/R,\sigma} \simeq \mathrm T_{A'/R',\sigma'}.
\end{displaymath}
\item If $B$ is a twisted localization of $A$, then there exists a canonical isomorphism
\begin{displaymath}
B \otimes_{A} \mathrm T_{A/R,\sigma} \simeq \mathrm T_{B/R,\sigma}. \quad \Box
\end{displaymath}
\end{enumerate}
\end{cor}

\begin{dfn} \label{twico} 
A \emph{twisted connection} on an $A$-module $M$ is an $R$-linear map
\begin{displaymath}
\nabla : M \mapsto M \otimes_{A} \Omega^1_{A/R,\sigma}
\end{displaymath}
such that
\begin{displaymath}
\forall s \in M, \forall x \in A, \quad \nabla (xs) = s \otimes \mathrm d(x) + \sigma(x)\nabla(s).
\end{displaymath}
An $A$-linear map between two $A$-modules with twisted connections is said to be \emph{horizontal} if it is compatible with the connections.
\end{dfn}

Clearly, $A$-modules endowed with a connection and horizontal maps form a category $\nabla_{\sigma_{A}}\mathrm{-Mod}$.

\begin{rmk}
This definition is compatible with definition 2.2.1 by Andr\'e in \cite{Andre01}.
In particular, all the tannakian formalism applies but this is not what we are interested in.
\end{rmk}

Recall from \cite{LeStumQuiros15a} that if $D$ is a twisted derivation of $A$, then a twisted $D$-derivation of an $A$-module $M$ is an $R$-linear endomorphism $D_{M}$ of $M$ that satisfies the {twisted Leibnitz rule}:
\begin{displaymath}
\forall x \in A, \forall s \in M, \quad D_{M}(xs) = D(x)s + \sigma_{A}(x)D_{M}(s).
\end{displaymath}
One may then consider the notion of \emph{action by twisted derivations} of $T_{A/R,\sigma}$ on $M$: it is an $R$-linear action such that whenever $D \in  T_{A/R,\sigma}$, the map $D_{M} : s \mapsto D.s$ is a $D$-derivation.

\begin{prop} \label{dact}
There exists a functor from the category of $A$-modules endowed with a twisted connection to the category of $A$-modules endowed with a linear action of $T_{A/R,\sigma}$ by twisted derivations.
It is an equivalence (an isomorphism) when $\Omega^1_{A/R,\sigma}$ is free of finite rank.
\end{prop}

\begin{pf}
If $M$ is endowed with a twisted connection $\nabla : M \mapsto M \otimes_{A} \Omega^1_{\sigma}$
 and $D$ is a twisted derivation of $A$, we may write uniquely $D = u \circ d$ with $u : \Omega^1_{\sigma} \to A$ and consider the composite map $D_{M} := (\mathrm{Id}_{M} \otimes u) \circ \nabla : M \to M$.
 Then, we will have
\begin{eqnarray*}
D_{M}(xs) 
&=& (\mathrm{Id}_{M} \otimes u)(\nabla(xs)) = (\mathrm{Id}_{M} \otimes u)(s \otimes \mathrm d(x) + \sigma(x)\nabla(s))
\\
&=& (\mathrm u \circ d)(x)s +\sigma(x) (\mathrm{Id}_{M} \otimes u)(\nabla(s)) = D(x)s +\sigma(x) D_{M}(s).
\end{eqnarray*}
Conversely, assume that $M$ is endowed with an action of $T_{\sigma}$ by twisted derivations.
Let $D_{1}, \ldots, D_{n}$ be a basis of $\mathrm{T}_{\sigma}$ and $\omega_{1}, \ldots, \omega_{n}$ be the dual basis in $\Omega^1_{\sigma}$.
Then, we can define $\nabla(s) = \sum D_{i,M}(s) \otimes \omega_{i}$ and check that
\begin{eqnarray*}
\nabla (xs)
&=& \sum D_{i,M}(xs) \otimes \omega_{i} = \left(\sum D_{i}(x)s +\sigma(x) D_{i,M}(s)\right) \otimes \omega_{i}
\\
&=& s \otimes  \sum D_{i}(x)\omega_{i} +\sigma(x)  \sum D_{i,M}(s) \otimes \omega_{i} = s \otimes \mathrm d(x) + \sigma(x)\nabla(s).
\end{eqnarray*}

Clearly, this is an inverse to the previous functor.
$\quad \Box$
\end{pf}

\begin{prop} \label{1coor} 
If $x \in A$ is a \emph{twisted coordinate} on $A$, then $\Omega^1_{A/R,\sigma}$ is free of rank $1$ generated by $\mathrm dx$.
Moreover, there exists a unique twisted derivation $\partial_{x,\sigma}$ of $A$ such that $\partial_{x,\sigma}(x) =1$ and we have
\begin{displaymath}
\forall D \in T_{A/R, \sigma}, \quad D = D(x)\partial_{x,\sigma}.
\end{displaymath}
\end{prop}

\begin{pf}
The first assertion is a particular case of corollary \ref{free}.
The second one then follows from corollary \ref{dualid}.
$\quad \Box$
\end{pf}

In particular, we see that a twisted coordinate is also a coordinate in the sense of \cite{LeStumQuiros15a}.
In order to lighten the notations, we will usually drop the index $x$ but we must not forget that $\partial_{\sigma}$ depends on the choice of $x$.
Also, we would rather write $\partial_{A,\sigma}$ than $\partial_{\sigma_{A}}$ when we want to make clear the dependence on $A$.

If $x$ is a twisted coordinate on $A$, one may consider the \emph{twisted Weyl algebra} $\mathrm D_{A/R,\sigma,\partial}$ (see \cite{LeStumQuiros15a} for example), that we will usually denote by $\mathrm D_{A/R,\sigma}$ and sometimes simply by $\mathrm D_{\sigma}$.
This is the non commutative polynomial ring in one variable $\partial_{\sigma}$ over $A$ with the commutation rule
\begin{displaymath}
\forall z \in A, \quad \partial_{\sigma} z = \partial_{A,\sigma}(z) + \sigma_{A}(z)\partial_{\sigma}.
\end{displaymath}
Moreover, there exists an equivalence (an isomorphism) of categories $\mathrm D_{A/R,\sigma}\mathrm{-Mod} \simeq \partial_{A,\sigma}\mathrm{-Mod}$ where the later denotes the category of $A$-modules 
$M$ endowed with a $\partial_{A,\sigma}$-derivation.
\begin{prop} \label{eqnabd}
Assume that $x$ is a twisted coordinate on $A$.
Then, there exists an equivalence (an isomorphism) of categories
\begin{displaymath}
\nabla_{\sigma_{A}}\mathrm{-Mod} \simeq D_{A/R,\sigma}\mathrm{-Mod}.
\end{displaymath}
\end{prop}

\begin{pf}
Follows from proposition \ref{dact}.
$\quad \Box$
\end{pf}

\section{Twisted binomial coefficient theorem for principal parts}

We prove here the main theorem that will allow us to recover twisted differential operators from principal parts. 
We use same notations as before.

In section \ref{tpp} we introduced the canonical map (it is a morphism of $R$-algebras)
\begin{displaymath}
\xymatrix@R0cm{
P_{A/R} \ar[r]^-{\delta_{A/R}} &P_{A/R} \otimes_{A} P_{A/R}
\\
x \tilde y \ar@{|->}[r] & x \otimes \tilde y
}.
\end{displaymath}
We want to investigate the interaction between $\sigma_{A}$ and $\delta_{A/R}$.

Recall that we also considered in section \ref{tpp} the maps
\begin{displaymath}
\xymatrix@R0cm{
P_{A/R} \ar[r]^-{\sigma_{P}} &P_{A/R}
\\
x \tilde y \ar@{|->}[r] & \sigma_{A}(x) \tilde y
},
\end{displaymath}
which is an $R$-linear ring homomorphism, and
\begin{displaymath}
\xymatrix@R0cm{
A \ar[r]^-{d_{A/R}} & P_{A/R}
\\
x \ar@{|->}[r] & \tilde x - x
},
\end{displaymath}
which is only $R$-linear.
As usual, we will drop the subscripts in order to lighten the notations, hoping that the meaning will always be clear from the context.

\begin{lem} \label{decompg}
For all $i = 0, \ldots, n$, we have in $P \otimes_{A} P$:
\begin{displaymath}
\forall x \in A, \quad \delta(\sigma^n(d(x))) = 1 \otimes \sigma^i(d(x)) + \sigma^{n-i}(d(\sigma^i(x)))\otimes 1.
\end{displaymath}
\end{lem}

\begin{pf}
We do the computations in $A \otimes_{R} A \otimes_{R} A$.
The right hand side is
\begin{align*}
&1 \otimes (1 \otimes x - \sigma^i(x) \otimes 1) + \sigma^{n-i}(1 \otimes \sigma^i(x) - \sigma^i(x) \otimes 1) \otimes 1
\\
&=1 \otimes 1 \otimes x - 1 \otimes \sigma^i(x) \otimes 1 + 1 \otimes \sigma^i(x) \otimes 1 - \sigma^n(x) \otimes 1 \otimes 1
\\
&= 1 \otimes 1 \otimes x - \sigma^n(x) \otimes 1 \otimes 1.
\end{align*}
If we develop the left hand side, we obtain exactly the same thing:
\begin{displaymath}
\delta(\sigma^n(d(x))) = \delta(1 \otimes x - \sigma^n(x) \otimes 1) = 1 \otimes 1 \otimes x - \sigma^n(x) \otimes 1 \otimes 1. \quad \Box
\end{displaymath}
\end{pf}

We endow $P \otimes_{A} P$ with the endomorphism $\sigma_{P} \otimes_{A} \mathrm{Id}_{P}$ (which is the same thing as $\sigma_{A} \otimes_{R} \mathrm{Id}_{A} \otimes_{R} \mathrm{Id}_{A}$ on $A \otimes_{R} A \otimes_{R} A$).

\begin{prop} \label{sigdel}
The map $\delta : P \to P \otimes_{A} P$ is a morphism of twisted $R$-algebras.
\end{prop}

\begin{pf}
This is the case $n = 1$ and $i = 0$ of lemma \ref{decompg}.
More precisely, if $x \in A$ and $\xi = \tilde x - x$, we have
\begin{displaymath}
\delta(\sigma(\xi)) =  1 \otimes \xi + \sigma(\xi) \otimes  1 = \sigma(\delta(\xi)). \quad \Box
\end{displaymath}
\end{pf}

\begin{prop} We have in $P \otimes_{A} P$:
\begin{displaymath}
\forall n \in \mathbb N, \quad \delta(I^{(n)}) \subset \sum_{i=0}^{n} I^{(i)} \otimes I^{(n-i)}.
\end{displaymath}
\end{prop}

\begin{pf}
First of all, since $I$ is generated by the image of $d$, it follows from lemma \ref{decompg} that for all $i = 0, \ldots, n$, we have
\begin{displaymath}
\delta(\sigma^n(I)) \subset P \otimes \sigma^i(I) + \sigma^{n-i}(I) \otimes P.
\end{displaymath}
Using induction, we obtain
\begin{eqnarray*}
\delta(I^{(n+1)})
&=& \delta(I^{(n)})\delta(\sigma^{n}(I))) \subset \sum_{i=0}^{n} (I^{(i)} \otimes I^{(n-i)})(P \otimes \sigma^{n-i}(I) + \sigma^i(I) \otimes P)
\\
&\subset& \sum_{i=0}^{n} (I^{(i)} \otimes I^{(n-i+1)}) + \sum_{i=0}^{n} (I^{(i+1)} \otimes I^{(n-i)}) \subset \sum_{i=0}^{n+1} I^{(i)} \otimes I^{(n+1-i)}. \quad \Box
\end{eqnarray*}
\end{pf}

\begin{cor} \label{delstab}
For all $m, n \in \mathbb N$, we have in $P \otimes_{A} P$:
\begin{displaymath}
\delta(I^{(n+m+1)}) \subset P \otimes_{A} I^{(m+1)} + I^{(n+1)} \otimes_{A} P.
\end{displaymath}
In other words, $\delta$ induces a map
\begin{displaymath}
\xymatrix@R0cm{
P_{(n+m)_{\sigma}} \ar[r]^-{\delta_{n,m}} &P_{(n)_{\sigma}} \otimes_{A} P_{(m)_{\sigma}}.
}
\end{displaymath}
\end{cor}

\textbf{Proof}
If $0 \leq i \leq m+n+1$, we have either $i > n$ and then $I^{(i)} \subset I^{(n+1)}$ or else $i \leq n$ so that $m + n + 1 -i > m$ and then $I^{(m+n+1-i)} \subset I^{(m+1)}$.
$\quad \Box$

Going to the limit, we obtain a canonical homomorphism of $R$-algebras
\begin{displaymath}
\xymatrix@R0cm{
\widehat P_{\sigma} \ar[r]^-{\widehat \delta} & \widehat P_{\sigma} \widehat \otimes_{A} \widehat P_{\sigma}
},
\end{displaymath}
where the right hand side is, \emph{by definition}, the inverse limit of all the $P_{(n)_{\sigma}} \otimes_{A} P_{(m)_{\sigma}}$.
In other words, we obtain a comultiplication on $\widehat P_{\sigma}$ that will allow us to turn its ``dual'' into a ring (more on this later).

We finish this section with the \emph{quantum binomial theorem for principal parts}:

\begin{thm} \label{biconf}
let $A$ be a twisted commutative $R$-algebra and $x \in A$ such that $\sigma(x) = qx + h$ with $q, h \in R$.
If we set $\xi = \tilde x - x$, then, we have
\begin{displaymath}
\delta(\xi^{(n)}) := \sum_{i=0}^{n} {n \choose i}_{q} \xi^{(n-i)} \otimes \xi^{(i)}.
\end{displaymath}
\end{thm}

\begin{pf}
The formula is proved to be correct by induction on $n$.
First of all, since $\delta$ is a ring homomorphism, we have
\begin{displaymath}
\delta(\xi^{(n+1)}) = \delta(\xi^{(n)}\sigma^{n}(\xi)) = \delta(\xi^{(n)})\delta(\sigma^{n}(\xi))).
\end{displaymath}
Using induction and lemma \ref{decomp} below, we get
\begin{eqnarray*}
\delta(\xi^{(n+1)})
&=& \sum_{i=0}^{n} {n \choose i}_{q} (\xi^{(n-i)} \otimes \xi^{(i)})(1 \otimes \sigma^i(\xi) + q^i\sigma^{n-i}(\xi)\otimes 1)
\\
&=& \sum_{i=0}^{n} {n \choose i}_{q} \xi^{(n-i)} \otimes \xi^{(i)}\sigma^{i}(\xi) + \sum_{i=0}^{n} {n \choose i}_{q} q^i\xi^{(n-i)}\sigma^{n-i}(\xi)  \otimes \xi^{(i)}
\\
&=& \sum_{i=0}^{n} {n \choose i}_{q} \xi^{(n-i)} \otimes \xi^{(i+1)} + \sum_{i=0}^{n} {n \choose i}_{q} q^i\xi^{(n-i+1)}  \otimes \xi^{(i)}
\\
&=& \sum_{i=1}^{n+1} {n \choose i-1}_{q} \xi^{(n-i+1)} \otimes \xi^{(i)} + \sum_{i=0}^{n} {n \choose i}_{q} q^i\xi^{(n-i+1)}  \otimes \xi^{(i)}
\\
&=& \sum_{i=0}^{n+1} \left({n \choose i-1}_{q} + q^i  {n \choose i}_{q}\right) \xi^{(n+1-i)} \otimes \xi^{(i)}
\\
&=& \sum_{i=0}^{n+1} {n +1 \choose i}_{q} \xi^{(n+1-i)} \otimes \xi^{(i)}. \quad \Box
\end{eqnarray*}
\end{pf}

\begin{lem} \label{decomp}
Under the hypothesis of the proposition, we have for all $i = 0, \ldots n$,
\begin{displaymath}
\delta(\sigma^n(\xi)) = 1 \otimes \sigma^i(\xi) + q^i\sigma^{n-i}(\xi)\otimes 1.
\end{displaymath}
\end{lem}

\begin{pf}
We have
\begin{displaymath}
d(\sigma(x)) = \widetilde{\sigma(x)} - \sigma(x)
= \widetilde {qx+h} - (qx + h) = q (\tilde x - x) = q \xi.
\end{displaymath}

The analogous result holds for the endomorphism $\sigma^i$.
It follows that $d(\sigma^i(x)) = q^i \xi$ and we finish with lemma \ref{decompg}.
$\quad \Box$
\end{pf}

\section{Twisted differential operators of infinite level}\label{tdo}

We are now able to define the ring of twisted differential operators (of infinite level).
We keep the previous notations.

\begin{dfn} \label{trudef}
If $M$ and $N$ are two $A$-modules, then a \emph{twisted differential operator $\varphi : M \to N$ of order at most $n$ (and infinite level)} is an $R$-linear map whose $A$-linearization 
\begin{displaymath}
\xymatrix@R0cm{ P \otimes_{A} M \ar@/^.5cm/[rr] ^-{\tilde \varphi} \ar[r]^{\simeq} &A \otimes_{R} M \ar[r] & N
\\ & x \otimes s \ar@{|->}[r] & x\varphi (s)}
\end{displaymath}
factors through $P_{(n)_{\sigma}}  \otimes_{A} M$.
\end{dfn}

Note that the condition means that the restriction of $\tilde \varphi$ to $I^{(n+1)} \otimes_{A} M$ is zero.
We might still write $\tilde \varphi$ for the map induced on $P_{(n)_{\sigma}}  \otimes_{A} M$ when there is no risk of confusion.

We denote by $\mathrm{Diff}_{n,\sigma}(M, N)$ the set of all twisted differential operators of order at most $n$.
Thus, we have a canonical isomorphism
\begin{displaymath}
\mathrm{Diff}_{n,\sigma}(M, N) \simeq \mathrm{Hom}_{A}(P_{(n)_{\sigma}} \otimes_{A} M, N).
\end{displaymath}
where $P_{(n)_{\sigma}}$ is seen as an $A$-module for the action on the \emph{right} with respect to $\otimes_{A}$ and for the action on the \emph{left} with respect to $\mathrm{Hom}_{A}$.
In particular, $\mathrm{Diff}_{n,\sigma}(M, N)$ has a natural structure of $P_{(n)_{\sigma}}$-module given by $x\tilde y \cdot \varphi := x \circ \varphi \circ y$, where multiplication by $y$ takes place in $M$ while we multiply by $x$ in $N$.

We will also denote by $\mathrm{Diff}_{\sigma}(M, N)$ the set of all twisted differential operators of any order so that
\begin{displaymath}
\mathrm{Diff}_{\sigma}(M, N) \simeq \varinjlim \mathrm{Hom}_{A}(P_{(n)_{\sigma}} \otimes_{A} M, N).
\end{displaymath}
In particular, we see that $\mathrm{Diff}_{\sigma}(M, N)$ has a natural structure of $\widehat P_{\sigma}$-module.

In the case $N=M$, we will write $\mathrm{Diff}_{n,\sigma}(M)$ and $\mathrm{Diff}_{\sigma}(M)$.
One also sets $\mathrm D^{(\infty)}_{A/R,\sigma} := \mathrm{Diff}_{\sigma}(A)$ and we will often drop the index $A/R$ and simply write $\mathrm D^{(\infty)}_{\sigma}$.

\begin{dfn} \label{deka}
Let $x$ be a twisted coordinate on $A$ and $\xi = \tilde x - x$.
Then, the \emph{standard basis} of $\mathrm{Diff}_{n,\sigma}(A)$ is the basis $\partial^{[k]}_{\sigma}$ dual to the images of the $\xi^{(k)}$ in $P_{(n)_{\sigma}}$.
We call $\partial^{[k]}_{\sigma}$ the \emph{standard twisted divided differential operator of order $k$ associated to $x$}.
\end{dfn}

Thus, when $x$ is a twisted coordinate, any $\varphi \in D_{\sigma}^{(\infty)}$ can be uniquely written as a finite sum $\sum z_{k} \partial_{\sigma}^{[k]}$ with $z_{k} \in A$ (and conversely, any such sum is in $D_{\sigma}^{(\infty)}$).
The canonical basis is characterized by the property
\begin{displaymath}
\forall k, n \in \mathbb N, \quad \widetilde{\partial^{[k]}_{\sigma}}(\xi^{(n)}) = \left\{\begin{array} {cll} 1 & \mathrm{if} \ k = n \\ 0 & \mathrm{otherwise}\end{array} \right. .
\end{displaymath}

We will give explicit examples later on.

The next proposition shows that the $A$-module of twisted differential operators could have been also defined by induction on the order $n$ (this is sometimes more convenient and does not require to work out the theory of principal parts).
Instead, it uses the notion of twisted bracket
\begin{displaymath}
\forall \varphi \in \mathrm{Hom}_{R}(M, N), \forall x \in A, \quad [\varphi, x]_{\sigma} = \varphi \circ x - \sigma(x) \circ \varphi
\end{displaymath}
already used in \cite{LeStumQuiros15a}.
We will need this intermediate result:

\begin{lem} \label{phipsi}
Let $M,N$ be two $A$-modules, $\varphi \in \mathrm{Hom}_{R}(M,N)$ and $x \in A$.
If we set $\varphi_{x} := [\varphi, x]_{\sigma}$,  then,
\begin{displaymath}
\tilde \varphi_{x} = \tilde \varphi \circ (\sigma(\xi) \otimes_{A} \mathrm{Id}_{M}) : P \otimes_{A} M \to N.
\end{displaymath}
\end{lem}

\begin{pf} We do the computations in $A \otimes_{R} M$.
Let $y \in A$ and $s \in M$, we have
\begin{eqnarray*}
\tilde \varphi (\sigma(\xi)(y \otimes s))
&=& \tilde \varphi ((1 \otimes x - \sigma(x) \otimes 1)(y \otimes s))= \tilde \varphi (y \otimes xs) - \tilde \varphi (\sigma(x)y \otimes s)
\\
&=& y \varphi(xs) - \sigma(x)y\varphi(s) = y \varphi_{x}(s) = \tilde \varphi_{x}(y \otimes s).
\quad \Box
\end{eqnarray*}
\end{pf}

\begin{prop}\label{carac}
Let $M,N$ be two $A$-modules and $\varphi \in \mathrm{Hom}_{R}(M,N)$.
Then, we have
\begin{displaymath}
\forall n \in \mathbb N, \quad \varphi \in \mathrm{Diff}_{n,\sigma}(M, N) \Leftrightarrow \forall x \in A, [\varphi,x]_{\sigma^n} \in \mathrm{Diff}_{n-1,\sigma}(M, N).
\end{displaymath}
\end{prop}

Note that we can start the induction process with $\mathrm{Diff}_{0,\sigma}(M, N)= \mathrm{Hom}_{A}(M, N)$ or with $\mathrm{Diff}_{-1,\sigma}(M, N) = 0$ if we prefer.

\begin{pf}
For $x \in A$, we set $\varphi_{x} := [\varphi, x]_{\sigma^n}$.
Then, we consider the linearizations
\begin{displaymath}
\tilde{\varphi}, \tilde \varphi_{x} : P \otimes_{A} M \to N
\end{displaymath}
of $\varphi$ and $\varphi_{x}$, respectively, and apply lemma \ref{phipsi} to $\sigma^n$ so that 
\begin{displaymath}
\tilde \varphi_{x} = \tilde \varphi \circ (\sigma^n(\xi) \otimes_{A} \mathrm{Id}_{M}).
\end{displaymath}
Thus we see that $\tilde \varphi = 0$ on $I^{(n+1)} \otimes_{A} M = I^{(n)}\sigma^n(I) \otimes_{A} M$ if and only if $\tilde \varphi_{x} = 0$ on $I^{(n)} \otimes_{A} M$ for all $x \in A$.
In other words, $\tilde \varphi$ factors through $P_{(n)_{\sigma}} \otimes_{A} M$ if and only if all $\tilde \varphi_{x}$ factor through $P_{(n-1)_{\sigma}} \otimes_{A} M$.
$\quad \Box$
\end{pf}

\begin{cor} \label{sigmadif}
A twisted differential operator of order at most $n$ (and infinite level) from $M$ to $N$ is an $R$-linear map $\varphi : M \to N$ such that
\begin{displaymath}
\forall x_{0}, \ldots, x_{n}\in A, \quad\ [[\cdots[ [\varphi,x_{n}]_{\sigma^n}, x_{n-1}]_{\sigma^{n-1}} \cdots ]_{\sigma}, x_{0}] = 0. \quad \Box
\end{displaymath}
\end{cor}

\begin{rmks}
\begin{enumerate}
\item
Be careful that, with the notations of Valery Lunts and Alexander L. Rosenberg in \cite {LuntsRosenberg97}, our $\mathrm{Diff}_{\sigma}(M, N)$ is different from their $\mathrm{Diff}(M, N^\sigma)$ which is defined by the condition
\begin{displaymath}
\forall x_{0}, \ldots, x_{n}\in A, \quad [ [\cdots[ [\varphi,x_{n}]_{\sigma}, x_{n-1}]_{\sigma} \cdots ]_{\sigma}, x_{0}] = 0.
\end{displaymath}
They only coincide when $n =0, 1$.
\item Our $\mathrm{Diff}_{\sigma}(M, N)$ should however coincide with some flavor of the $\mathrm{Diff}_{\beta}(M, N)$ of Lunts and Rosenberg.
More precisely, in order to define this module, they need a $G$-grading on $A$ and a bilinear map $\beta : G \times G \to R^\times$ (they use $k$ and $R$ for our $R$ and $A$).
In the simplest non trivial case $A = R[x]$ and $\beta(m, n) = q^{-mn}$, we believe that their $D_{\beta}(A)$ coincides with our $D_{A/R,\sigma}^{(\infty)}$ but their $D_{q}(A)$ is bigger (see \cite{IyerMcCune02} for example).
\item Charlotte Hardouin introduces in \cite{Hardouin10}, definition 2.4, what she calls an \emph{iterative $q$-difference ring} or $\mathrm{ID}_{q}$-ring for short.
She chooses some non zero $q \in K$ where $K$ is a fixed algebraically closed field and endows the field $A := K(x)$ of rational functions on $K$ with the automorphism $\sigma(x) = qx$.
Then, with our notations, an $\mathrm{ID}_{q}$-ring is a finitely generated $A$-algebra $B$ with a structure of $D_{A/K,\sigma}^{(\infty)}$-module, denoted by $(\varphi, y) \mapsto \varphi(y)$, such that the map $y \mapsto \sigma (y)$ is an automorphism of the ring $B$ and
\begin{equation} \label{hardouin}
\forall k \in \mathbb N, \forall y, z \in B, \quad \partial^{[k]}_{\sigma}(yz) = \sum_{i+j=k} (\sigma^i\partial^{[j]}_{\sigma})(y)\partial^{[i]}_{\sigma}(z).
\end{equation}
Note that $B$ becomes an inversive twisted $A$-algebra and that condition \eqref{hardouin} is automatic if the $q$-characteristic of $B$ is zero.
\end{enumerate}
\end{rmks}

\begin{prop} \label{compo}
Composition of twisted differential operators gives a twisted differential operator.
Moreover, its order is at most the sum of the order of the components.
\end{prop}

\begin{pf}
We let $\varphi : M \to N$ be a twisted differential operator of order $n$ and $\psi : L \to M$ a twisted differential operator of order $m$ and consider the factorization
\begin{equation*}
\xymatrix@R0cm{ P \otimes_{A} L \ar[r]^-{\delta \otimes \mathrm{Id}} & P \otimes_{A} P \otimes_A L \ar[r]^-{\mathrm{Id} \otimes_{A} \tilde \psi} & P \otimes_{A} M \ar[r]^-{\tilde \varphi} & N}
\end{equation*}
of lemma \ref{compolin}.
The map $\tilde \varphi$ factors through $P_{(n)_{\sigma}} \otimes_{A} M$ and $\mathrm{Id} \otimes \tilde \psi$ factors through $P \otimes_{A} P_{(m)_{\sigma}} \otimes_{A} L$. Therefore, their composite factors through $P_{(n)_{\sigma}} \otimes_{A} P_{(m)_{\sigma}} \otimes_{A} L$ and it follows from corollary \ref{delstab} that the whole thing will factor through $P_{(n+ m)_{\sigma}} \otimes L$.
Thus, $\varphi \circ \psi$ is a twisted differential operator of order at most $m + n$.
$\quad \Box$
\end{pf}

\begin{prop}
\begin{enumerate}
\item
Let $R \to R'$ be any morphism of commutative rings and $A' := R' \otimes_R A$, endowed with $\mathrm{Id}_{R'} \otimes_{R} \sigma_{A}$.
Then, we have
\begin{displaymath}
\mathrm{Diff}_{\sigma_{A'}}(A' \otimes_A M, A' \otimes_A N) \simeq A' \otimes_A \mathrm{Diff}_{\sigma_{A}}(M, N).
\end{displaymath}
\item  If $A \to B$ is a localization of twisted $R$-algebras, we have
\begin{displaymath}
\mathrm{Diff}_{\sigma_{B}}(B \otimes_A M, B \otimes_A N) \simeq B \otimes_A \mathrm{Diff}_{\sigma_{A}}(M, N).
\end{displaymath}
\end{enumerate}
\end{prop}

\begin{pf}
Follows from propositions \ref{extn} and \ref{locord} and the fact that direct limits commute with tensor product.
$\quad \Box$
\end{pf}

\begin{rmk}
\begin{enumerate}
\item As a particular case, we will have
\begin{displaymath}
\mathrm D^{(\infty)}_{A'/R',\sigma'} \simeq A' \otimes_A \mathrm D^{(\infty)}_{A/R,\sigma}\quad \mathrm{and} \quad \mathrm D^{(\infty)}_{B/R,\sigma} \simeq B \otimes_A \mathrm D^{(\infty)}_{A/R,\sigma}.
\end{displaymath}
\item When $A$ is a quantum $R$-algebra, we will always have (see the remark following proposition \ref{locco}) a natural isomorphism of $A$-modules
\begin{displaymath}
\mathrm D^{(\infty)}_{A/R,\sigma} \simeq A \otimes_{R[x]} \mathrm D^{(\infty)}_{R[x]/R,\sigma}.
\end{displaymath}
Note however that the ring structure (or equivalently the action on $A$) plays a fundamental role.
\end{enumerate}
\end{rmk}

\begin{prop}
We have $\mathrm{Diff}_{0,\sigma}(A) = A$ and $\mathrm{Diff}_{1,\sigma}(A) = A \oplus \mathrm T_{A/R,\sigma}$.
\end{prop}

\begin{pf}
The first assertion follows from the fact that $P_{(0)_{\sigma}} = A$ and the second one from proposition \ref{splitex}.
$\quad \Box$
\end{pf}

Recall that we introduced in \cite{LeStumQuiros15a} the ring $\overline{\mathrm D}_{A/R,\sigma}$ of \emph{small} twisted differential operators of $A/R$ as the smallest subring of $\mathrm{End}_{R}(A)$ containing both $A$ and $\mathrm T_{A/R,\sigma}$.
Again, we will simply write $\overline{\mathrm D}_{\sigma}$ when we believe that there is no risk of confusion.
Then, we have the following:

\begin{cor} \label{contin}
The ring of small twisted differential operators is contained inside the ring of twisted differential operators of infinite level:
we have
\begin{displaymath}
\overline{\mathrm D}_{A/R,\sigma} \subset \mathrm D^{(\infty)}_{A/R,\sigma}. \quad \Box
\end{displaymath}
\end{cor}

When there exists a twisted coordinate, we can make the twisted Weyl algebra enter the picture and we have:

\begin{cor} \label{epimon}
If $x$ is a twisted coordinate on $A$, there exists an epi-mono factorization
\begin{displaymath}
\xymatrix@R0cm{ \mathrm D_{A/R,\sigma} \ar@{->>}[r] & \overline{\mathrm D}_{A/R,\sigma} \ar@{^{(}->}[r] & \mathrm D^{(\infty)}_{A/R,\sigma}.
}
\end{displaymath}
\end{cor}

\begin{pf}
There exists a natural map $\mathrm D_{\sigma} \to \overline{\mathrm D}_{\sigma}$ that sends the parameter $\partial_{\sigma}$ of $\mathrm D_{\sigma,\partial}$ to the corresponding endomorphism $\partial_{A,\sigma}$ of $A$.
And it is surjective since $\partial_{A,\sigma}$ is a generator of $\mathrm{T}_{\sigma}$.
$\quad \Box$
\end{pf}

At some point, we will need to be able to compare twisted differential operators with respect to $\sigma$ and twisted differential operators with respect to the powers (or roots) of $\sigma$.

\begin{prop} \label{pow}
For all $m > 0$, if $M, N$ are two $A$-modules, we have
\begin{displaymath}
\mathrm{Diff}_{\sigma^m}(M, N) \subset \mathrm{Diff}_{\sigma}(M, N).
\end{displaymath}
\end{prop}

\begin{pf}
Since
\begin{displaymath}
I^{(mn)_{\sigma}} = I \sigma(I) \cdots \sigma^{mn}(I) \subset I\sigma^m(I) \cdots \sigma^{mn}(I) = I^{(n)_{\sigma^m}},
\end{displaymath}
there exists an natural surjective map
\begin{displaymath}
P_{(nm)_{\sigma}} \to P_{(n)_{\sigma^m}},
\end{displaymath}
from which we derive an inclusion $\mathrm{Diff}_{n, \sigma^m}(M, N) \subset \mathrm{Diff}_{mn,\sigma}(M, N)$.
$\quad \Box$
\end{pf}

Recall from \cite{LeStumQuiros15} that a system of roots of $\sigma$ is a family $\underline \sigma := \{\sigma_{n}\}_{n \in S}$, with $S \subset \mathbb N$, of $R$-linear ring endomorphisms of $A$ such that $\sigma_{n}^m = \sigma_{n'}^{m'}$ whenever $m/n = m'/n'$ and $\sigma_{n}^n = \sigma$.
We will always assume that $S$ is filtering for division.
We will call the pair $(A, \underline \sigma)$ an \emph{$S$-twisted $R$-algebra}.

\begin{dfn}
Let $\underline \sigma := \{\sigma_{i}\}_{i \in S}$ be a system of roots of $\sigma$.
Then the \emph{ring of twisted differential operators (of infinite level)} $\mathrm D^{(\infty)}_{A/R,\underline \sigma}$ is the $R$-subalgebra of $\mathrm{End}_{R}(A)$ generated by all $\sigma_{n}$-differential operators (of infinite level) for all $n \in S$.
\end{dfn}

Proposition \ref{pow} then has the following consequence:

\begin{cor}
If $\underline \sigma$ is a system of roots of $\sigma$, we have
\begin{displaymath}
\mathrm D^{(\infty)}_{A/R,\underline \sigma} =\cup  \mathrm D^{(\infty)}_{A/R,\sigma_{n}}. \quad \Box
\end{displaymath}
\end{cor}

Note that in section 3 of \cite{LeStumQuiros15a} we defined exactly in the same way the ring of small twisted differential operators $\overline {\mathrm D}_{A/R,\underline \sigma}$ for any family $\underline \sigma$ (not necessarily a root system), but we showed that the analogous statement is \emph{not} true in general.

\section{Twisted taylor series}

We will develop here the formalism of twisted Taylor maps which describes the formal solutions of twisted differential modules.
Notations are as before.

\begin{lem} \label{disrep}
If $M$ is a $\mathrm D_{A/R,\sigma}^{(\infty)}$-module, then the canonical map $\mathrm D_{A/R,\sigma}^{(\infty)} \to \mathrm{End}_{R}(M)$ induces, for all $n \in \mathbb N$, a $P_{A/R,(n)_{\sigma}}$-linear map
\begin{displaymath}
\mathrm{Diff}_{n,\sigma}(A) \to \mathrm{Diff}_{n,\sigma}(M).
\end{displaymath}
Hence, there exists a canonical $\widehat P_{\sigma}$-linear map
\begin{displaymath}
D^{(\infty)}_{A/R} \to \mathrm{Diff}_{\sigma}(M).
\end{displaymath}
\end{lem}

\begin{pf}
Since the canonical map $\lambda : D_{\sigma}^{(\infty)} \to \mathrm{End}_{R}(M)$ is a morphism of $A$-algebras, it will commute with the action of $P$.
More precisely, for all $x,y \in A$ and $\varphi \in D_{\sigma}^{(\infty)}$, we have (see section \ref{tdo} for the definition of the action of $P$ on $D_{\sigma}^{(\infty)}$)
\begin{displaymath}
\lambda(x\tilde y \cdot \varphi) = \lambda(x \circ \varphi \circ y) = x \circ \lambda(\varphi) \circ y = x\tilde y \cdot \lambda(\varphi).
\end{displaymath}
In particular, if $\tilde \varphi$ is zero on $I^{(n+1)}$, then $\widetilde {\lambda(\varphi)}$ will be zero on $I^{(n+1)} \otimes_{A} M$.
It means that the image of $\mathrm{Diff}_{n,\sigma}(A)$ falls inside $\mathrm{Diff}_{n,\sigma}
(M)$.
$\quad \Box$
\end{pf}

We will usually denote by $\varphi_{M} \in \mathrm{Diff}_{\sigma}(M)$ the image of $\varphi \in D^{(\infty)}_{\sigma}$.
In other words, for $\varphi \in D_{\sigma}^{(\infty)}$ and $s \in M$, we will have $\varphi_{M}(s) = \varphi s$.

\begin{dfn}
A \emph{twisted Taylor structure (of infinite level)} on an $A$-module $M$ is a compatible family of $A$-linear maps $\theta_{n} : M \to M \otimes_{A} P_{(n)_{\sigma}}$ (called \emph{twisted Taylor maps}) with $\theta_{0} = \mathrm{Id}$, making commutative all the diagrams
\begin{displaymath}
\xymatrix{M \ar[rr]^{\theta_{m}} \ar[d]^{\theta_{n+m}} && M \otimes_{A} P_{(m)_{\sigma}} \ar[d]^{\theta_{n} \otimes \mathrm{Id}} \\ M \otimes_{A} P_{(n+m)_{\sigma}} \ar[rr]^{\mathrm{Id} \otimes \delta_{n,m}} && M \otimes_{A} P_{(n)_{\sigma}} \otimes_{A} P_{(m)_{\sigma}}}.
\end{displaymath}
\end{dfn}

For example, the canonical twisted Taylor structure on $A$ is defined by the family of composite maps
\begin{displaymath}
\xymatrix{ A \ar[r]^{\theta} \ar[dr]_{\theta_{n}} & P \ar@{->>}[d]
\\   & P_{(n)_{\sigma}}},
\end{displaymath}
where the upper map is the Taylor map $x\mapsto \tilde x$ given by the action on the right (see section \ref{tpp}).

There exists an obvious notion of morphism of $A$-modules endowed with a twisted Taylor structure and they form a category.

\begin{prop} \label{Taylorstruc}
Let $M$ be an $A$-module endowed with a twisted Taylor structure $(\theta_{n})_{n \in \mathbb N}$.
Then, there exists a unique structure of $\mathrm D^{(\infty)}_{A/R,\sigma}$-module on $M$ such that
\begin{equation} \label{Taylexp}
\forall n \in \mathbb N, \quad \theta_{n}(s) = \sum s_{k} \otimes f_{k} \Rightarrow \forall \varphi \in \mathrm{Diff}_{n,\sigma}(A), \varphi_{M}(s) = \sum \tilde \varphi(f_{k}) s_{k}.
\end{equation}
This is functorial in $M$ and we obtain an equivalence (an isomorphism) of categories if all $P_{(n)_{\sigma}}$ are finite projective (for the left $A$-module structure).
\end{prop}

Note that the last condition is satisfied if there exists a twisted coordinate on $A$.

\begin{pf}
First of all, there exists for all $n \in \mathbb N$, a canonical morphism of $A$-modules
\begin{equation}\label{dualm}
M \otimes_{A} P_{(n)_{\sigma}} \to \mathrm{Hom}_{A}(\mathrm{Hom}_{A}(P_{(n)_{\sigma}}, A), M).
\end{equation}
which is automatically $P_{(n)_{\sigma}}$-linear.
Now, $A$-linear maps
\begin{equation} \label{nd}
\theta_{n} : M \to M \otimes_{A} P_{(n)_{\sigma}}
\end{equation}
correspond bijectively to $P_{(n)_{\sigma}}$-linear maps
\begin{displaymath}
\epsilon_{n} : P_{(n)_{\sigma}} \otimes_{A} M \to M \otimes_{A} P_{(n)_{\sigma}},
\end{displaymath}
and we can compose with the map \eqref{dualm} in order to get
\begin{displaymath}
P_{(n)_{\sigma}} \otimes_{A} M \to \mathrm{Hom}_{A}(\mathrm{Hom}_{A}(P_{(n)_{\sigma}}, A), M),
\end{displaymath}
or equivalently,
\begin{displaymath}
\mathrm{Hom}_{A}(P_{(n)_{\sigma}}, A) \to \mathrm{Hom}_{A}(P_{(n)_{\sigma}} \otimes_{A} M, M).
\end{displaymath}
In other words, we obtain $P_{(n)_{\sigma}}$-linear maps
\begin{equation} \label{dndn}
\xymatrix@R0cm{ \mathrm{Diff}_{n,\sigma}(A) \ar[r] & \mathrm{Diff}_{n,\sigma}(M)
\\ \varphi \ar@{|->}[r] & \varphi_{M}}.
\end{equation}
Formula \eqref{Taylexp} follows directly from the construction.
Compatibility for various $n$ in \eqref{dndn} follows from compatibility for various $n$ in \eqref{nd}.
We need to show that the corresponding map $\mathrm D_{\sigma}^{(\infty)} \to \mathrm{End}(M)$ is a morphism of rings.
To do that, one can use the description of composition of twisted differential operators given in proposition \ref{compo}: we need to verify that the maps \eqref{nd} are compatible with $\delta$ which is exactly the condition in the definition of Taylor structure.

This construction is clearly functorial.
Moreover, if $P_{(n)_{\sigma}}$ is finite projective (for the left $A$-module structure), then the map \eqref{dualm} is actually an isomorphism.
And it follows from lemma \ref{disrep} that a $\mathrm D^{(\infty)}_{\sigma}$-module structure on $M$ will provide us with a compatible family of maps as in \eqref{dndn}.
$\quad \Box$
\end{pf}

\begin{dfn} \label{taylordf}
Let $M$ be an $A$-module endowed with a twisted Taylor structure.
Then, the \emph{twisted Taylor map} of $M$ is the map
\begin{displaymath}
\widehat \theta = \varprojlim \theta_{n} : M \mapsto M \widehat \otimes_{A} \widehat P_{\sigma} := \varprojlim  M \otimes_{A} P_{(n)_{\sigma}}.
\end{displaymath}
The \emph{twisted Taylor series} of $s \in M$ is $\widehat \theta(s) \in M \widehat \otimes_{A} \widehat P_{\sigma}$.
\end{dfn}

There exists a commutative diagram
\begin{displaymath}
\xymatrix{M \ar[rr]^{\widehat \theta} \ar[d]^{\widehat\theta} && M \widehat \otimes_{A} \widehat P_{\sigma} \ar[d]^{\widehat \theta \widehat \otimes \mathrm{Id}} \\ M \widehat \otimes_{A} \widehat P_{\sigma} \ar[rr]^{\mathrm{Id} \widehat \otimes \widehat \delta} && M \widehat \otimes_{A} \widehat P_{\sigma} \widehat \otimes_{A} \widehat P_{\sigma}},
\end{displaymath}
and the action of $A$ on $\widehat P_{\sigma}$ on the right is given by the Taylor map of $A$.

In practice, we will only consider the case of finitely presented $A$-modules $M$, and then, the completed tensor product is the usual tensor product $M \otimes_{A} \widehat P_{\sigma}$

We just showed in proposition \ref{Taylorstruc} that any $\mathrm D^{(\infty)}_{\sigma}$-module comes with a canonical twisted Taylor structure.
When there exists a twisted coordinate on $A$, we can can describe it explicitly as follows:

\begin{prop}[Twisted Taylor formula] \label{qtayl}
Assume that $x$ is a twisted coordinate on $A$ and let $\xi = \tilde x -x$.
If $M$ is a $\mathrm D^{(\infty)}_{\sigma}$-module, we have for all $n \in \mathbb N$,
\begin{displaymath}
\theta_{n}(s) = \sum_{k=0}^{n} \partial_{\sigma}^{[k]}(s) \otimes \xi^{(k)} \in  M  \otimes_{A} P_{(n)_{\sigma}}
\end{displaymath}
and
\begin{displaymath}
\widehat \theta(s) = \sum_{k=0}^{\infty} \partial_{\sigma}^{[k]}(s) \otimes \xi^{(k)} \in M \widehat \otimes_{A} \widehat P_{\sigma}.
\end{displaymath}
\end{prop}

\begin{pf}
This follows from equation \eqref{Taylexp}: if we write
\begin{displaymath}
\theta_{n}(s) = \sum_{k=0}^{n} s_{k} \otimes \xi^{(k)},
\end{displaymath}
we will have for all $l \in \mathbb N$,
\begin{displaymath}
\partial_{\sigma}^{[l]}(s) = \sum_{k=0}^{n} \tilde \partial_{\sigma}^{[l]}(\xi^{(k)}) s_{k} = s_{l}.\quad \Box
\end{displaymath}
\end{pf}

In particular, we see that, if $z \in A$, then the image in $\widehat P_{\sigma}$ of $\tilde z = 1 \otimes z \in P$ is the twisted Taylor series
\begin{displaymath}
\widehat \theta(z) = \sum_{k} \partial_{\sigma}^{[k]}(z)\xi^{(k)}.
\end{displaymath}
This explains why it is legitimate to call the map $z \mapsto \tilde z$ the Taylor map.

\begin{xmps}
\begin{enumerate}
\item If $x$ is a twisted coordinate on $A$, we have $\widehat \theta(x) = x + \xi$ and
\begin{displaymath}
\widehat \theta(x^2) = x^2 + (x + \sigma(x)) \xi + \xi^2.
\end{displaymath}
\item Assume that $x \in A^\times$ is an invertible twisted coordinate on $A$ and that $\sigma(x) = qx$ with $q \in R^\times$, then we have
\begin{displaymath}
\widehat \theta \left( \frac 1x \right) = \sum_{k=0}^{\infty} (-1)^k \frac  {\xi^{(k)}}{ q^{\frac {k(k+1)}2} x^{k+1}} = \frac 1x - \frac \xi{qx^2} + \frac {\xi^{(2)}}{q^3x^3} - \cdots.
\end{displaymath}
\end{enumerate}
\end{xmps}
 
\begin{rmk}
If $A$ is a twisted localization of $R[x]$, then there exists at most one $R$-algebra homomorphism $\widehat \theta : A \to \widehat P_{\sigma}$ such that $\widehat \theta(x) = \tilde x$.
It means that the twisted Taylor map
\begin{displaymath}
\xymatrix@R=0cm{
A \ar[r]^{\widehat \theta  } & \widehat P_{\sigma} \\
z \ar@{|->}[r] & \sum_{k=0}^{\infty} \partial_{\sigma}^{[k]}(z) \xi^{(k)}}
\end{displaymath}
is the \emph{unique} such map.
\end{rmk}

\section{Quantum differential operators}

In the quantum situation, we can be a lot more explicit as we shall see now.
Thus, we assume in this section that $A$ is a quantum $R$-algebra:
we are given a twisted coordinate $x$ such that $\sigma(x) = qx+h$ with $q,h \in R$.

\begin{prop}
We have
\begin{displaymath}
\forall k,l \in \mathbb N, \quad \partial_{\sigma}^{[k]} \circ \partial_{\sigma}^{[l]} = {k+l \choose l}_{q} \partial_{\sigma}^{[k+l]}.
\end{displaymath}
\end{prop}

\begin{pf}
It follows from lemma \ref{compolin} that
\begin{displaymath}
\widetilde {\partial_{\sigma}^{[k]} \circ \partial_{\sigma}^{[l]}} = \widetilde {\partial_{\sigma}^{[k]}} \circ (\mathrm{Id} \otimes \widetilde {\partial_{\sigma}^{[l]}}) \circ \delta.
\end{displaymath}
Thus, using theorem \ref{biconf}, we see that
\begin{eqnarray*}
(\widetilde {\partial_{\sigma}^{[k]} \circ \partial_{\sigma}^{[l]}})(\xi^{(n)})
&=& \widetilde {\partial_{\sigma}^{[k]}} \left( (\mathrm{Id} \otimes \widetilde {\partial_{\sigma}^{[l]}}) \left( \sum_{i=0}^{n} {n \choose i}_{q} \xi^{(n-i)} \otimes \xi^{(i)}\right)\right)
\\
&=& \sum_{i=0}^{n}  {n \choose i}_{q} \widetilde {\partial_{\sigma}^{[k]}} \left( \xi^{(n-i)} \otimes \widetilde {\partial_{\sigma}^{[l]}}(\xi^{(i)})\right)
\\
&=&  {n \choose l}_{q} \widetilde {\partial_{\sigma}^{[k]}} (\xi^{(n-l)})
=  \left\{ \begin{array}{cl} {k+l \choose l}_{q} & \mathrm{if} \ n = k + l \\ 0 & \mathrm{otherwise}. \end{array} \right. . \quad \Box
\end{eqnarray*}
\end{pf}

Recall from \cite{LeStumQuiros15} that if $m \in \mathbb N$, we write $(m)_{q} := {m \choose 1}_{q} = 1 + \cdots + q^{m-1}$ and we define by induction $(m)_{q}! := (m)_{q} (m-1)_{q}!$.

\begin{cor} \label{formd}
We have
\begin{displaymath}
\forall k \in \mathbb N, \forall z \in A,  \quad \partial_{\sigma}^{k}(z) = (k)_{q}!  \partial_{\sigma}^{[k]}(z).
\end{displaymath}
\end{cor}

\begin{pf}
We proceed by induction on $k$ and obtain
\begin{displaymath}
 \partial_{\sigma}^{k+1}(z) = \partial_{\sigma}^{k} (\partial_{\sigma}(z)) = (k)_{q}!  \partial_{\sigma}^{[k]} (\partial_{\sigma}(z)) = (k)_{q}!  (\partial_{\sigma}^{[k]} \partial_{\sigma})(z)
 \end{displaymath}
 \begin{displaymath}
 = (k)_{q}! {k+1 \choose 1}_{q}  \partial_{\sigma}^{[k+1]} (z) = (k+1)_{q}!  \partial_{\sigma}^{[k+1]}(z) . \quad \Box
\end{displaymath}
\end{pf}

The next result is important because it describes explicitly the relations between the different rings of twisted differential operators introduced so far.
Recall from proposition \ref{epimon} that there exists an epi-mono factorization
\begin{displaymath}
\xymatrix@R0cm{ \mathrm D_{A/R,\sigma} \ar@{->>}[r] & \overline{\mathrm D}_{A/R,\sigma} \ar@{^{(}->}[r] & \mathrm D^{(\infty)}_{A/R,\sigma}.
}
\end{displaymath}

Recall from \cite{LeStumQuiros15} that the ring $R$ is said to be \emph{$q$-flat} (resp. \emph{$q$-divisible}) if $(m)_{q}$ is always regular (resp. invertible) unless $(m)_{q} = 0$ and also that the \emph{$q$-characteristic} of $R$ is the smallest positive integer $p$ such that $(p)_{q} = 0$, if it exists, and $0$ otherwise.

\begin{thm} \label{thd}
Assume $R$ is $q$-divisible and let $A$ be a $q$-$R$-algebra.
Then,
\begin{enumerate}
\item \label{qrat} If $q\mathrm{-char}(A) = 0$, we have
\begin{displaymath}
\mathrm D_{A/R,\sigma} = \overline{\mathrm D}_{A/R,\sigma} = \mathrm D^{(\infty)}_{A/R,\sigma}.
\end{displaymath}
\item \label{qlitl} If $q\mathrm{-char}(A) = p >0$, we have
\begin{displaymath}
\mathrm D_{A/R,\sigma}/\partial_{\sigma}^p \simeq \overline{\mathrm D}_{A/R,\sigma} \simeq \mathrm D^{(\infty)}_{A/R,\sigma}/K_{\sigma}^{[p]},
\end{displaymath}
where $K_{\sigma}^{[p]}$ is the free $A$-module generated by all $\partial_{\sigma}^{[k]}$ for $k \geq p$.
\end{enumerate}
\end{thm}

Note that with some extra conventions, assertions \ref{qrat}) and \ref{qlitl}) could be turned into a single one.

\begin{pf}
We use the formula from corollary \ref{formd}.

In situation \ref{qrat}), then all $q$-integers are invertible in $A$ and the composite map sends the canonical basis $\{\partial_{\sigma}^{k}\}_{k \in \mathbb N}$ of $\mathrm D_{\sigma}$ bijectively onto a basis of $\mathrm D^{(\infty)}_{\sigma}$.
The first assertion follows.

In situation \ref{qlitl}), we have $(p)_{q} = 0$ in $A$ and the element $\partial_{\sigma}^{p} \in \mathrm D_{\sigma}$ is therefore sent to $0$.
But since $R$ is $q$-divisible, we have $(m)_{q} \in R^\times$ for $m < p$, and the composite map sends the family $\{\partial_{\sigma}^{k}\}_{k < p}$ of $\mathrm D_{\sigma}$ bijectively onto a basis of $\mathrm D^{(\infty)}_{\sigma}/K_{\sigma}^{[p]}$.
$\quad \Box$
\end{pf}

\begin{rmks}
\begin{enumerate}
\item 
The hypothesis in \ref{qrat}) is satisfied in the cases of classical differential equations, classical finite difference equations and classical $q$-difference equations.

More precisely, it is satisfied for example when
\begin{enumerate}
\item $q = 1$ and $R$ is a $\mathbb Q$-algebra, or
\item  $q$ is not equal to zero, not a root of $1$ and belongs to a subfield $K$ of $R$.
\end{enumerate}
The hypothesis in \ref{qlitl}) is satisfied for differential equations and finite difference equations in positive characteristic as well as in the classical quantum case.
More precisely, they are satisfied for example when
\begin{enumerate}
\item $q = 1$ and $R$ is an $\mathbb F_{p}$-algebra, or
\item  $q$ is a non trivial $p$-th root of $1$ ($p$ not necessarily prime) and belongs to a subfield $K$ of $R$, or
\item  $q$ is a non trivial $p$-th root of $1$ with $p$ prime (but $q$ does not necessarily belong to a subfield of $R$).
\end{enumerate}
\item Since both $\mathrm D_{A/R,\sigma}$ and $\mathrm D^{(\infty)}_{A/R,\sigma}$ commute with extensions of $R$ (although  $\overline{\mathrm D}_{A/R,\sigma}$ does not), and $\mathrm D_{A/R,\sigma}$ always commutes with extensions of $A$, we can sometimes (see the remark following proposition \ref{locco}) reduce questions to the generic case
\begin{equation} \label{genc}
R = \mathbb Q(t)[s], \quad A = R[x], \quad q = t \quad \mathrm{and} \quad h = s,
\end{equation}
and work as well over the later.
In this case, thanks to the theorem, we can identify $\mathrm D_{A/R,\sigma}$ with $\mathrm D^{(\infty)}_{A/R,\sigma}$.
\item The same proof shows that if $A$ is a twisted $R$-algebra which is only $q$-flat (but not necessarily $q$-divisible), then
\begin{enumerate}
\item \label{qrat2} If $q\mathrm{-char}(A) = 0$, then $\mathrm D_{A/R,\sigma} = \overline{\mathrm D}_{A/R,\sigma}$.
\item \label{qlitl2} If $q\mathrm{-char}(A) = p>  0$, then $\mathrm D_{A/R,\sigma}/\partial_{\sigma}^p \simeq \overline{\mathrm D}_{A/R,\sigma}$.
\end{enumerate}
\end{enumerate}
\end{rmks}

The end of this section will be devoted to giving explicit formulas.
They are usually quite formal to prove in the ring of twisted differential operators of infinite level and their analog in the twisted Weyl algebra is then easily obtained thanks to theorem \ref{thd}.

It is useful to have a general formula for the commutation of twisted differential operators with the twisted coordinate $x$:

\begin{prop} \label{commix}
We have
\begin{displaymath}
\forall k \in \mathbb N \setminus \{0\}, \quad \partial_{\sigma}^{[k]} \circ x = \sigma^k(x) \partial_{\sigma}^{[k]} +  \partial_{\sigma}^{[k-1]}.
\end{displaymath}
\end{prop}

\begin{pf}
We start with the following computation:
\begin{displaymath}
\xi^{(n+1)} = \xi^{(n)}\sigma^n(\xi) = \xi^{(n)}(\tilde x - \sigma^n(x)) = \tilde x \xi^{(n)} - \sigma^n(x)\xi^{(n)},
\end{displaymath}
from which we derive
\begin{displaymath}
\tilde x \xi^{(n)} = \xi^{(n+1)} + \sigma^n(x)\xi^{(n)}.
\end{displaymath}
We have
\begin{displaymath}
\widetilde{\partial_{\sigma}^{[k]} \circ x} = \tilde \partial_{\sigma}^{[k]} \circ (\mathrm{ Id} \otimes \widetilde{m_x}) \circ \delta. 
\end{displaymath}
In this formula, $\widetilde{m_x}$  denotes the linearization of multiplication by $x$, which is easily seen to be given by
\begin{displaymath}
\widetilde{m_x}(\xi^{(i)}) = \left\{ \begin{array}{l} x \quad \mathrm{if} \quad i = 0 \\ 0 \quad \mathrm{if} \quad i \neq 0 \end{array} \right. .
\end{displaymath}
Now, we know from theorem \ref{biconf} (the quantic binomial formula for principal parts) that
\begin{displaymath}
\delta(\xi^{(n)}) = \sum_{i=0}^{n} {n \choose i}_{q} \xi^{(n-i)} \otimes \xi^{(i)}.
\end{displaymath}
It follows that
\begin{displaymath}
(\mathrm{ Id} \otimes \widetilde{m_x}) (\delta (\xi^{(n)})) =  \tilde x \xi^{(n)}.
\end{displaymath}
Putting all these together, we obtain
\begin{eqnarray*}
\widetilde{(\partial_{\sigma}^{[k]} \circ x)} (\xi^{(n)}) 
&=&  \tilde \partial_{\sigma}^{[k]} (\tilde x \xi^{(n)})
\\
&=& \tilde \partial_{\sigma}^{[k]} (\xi^{(n+1)} + \sigma^n(x)\xi^{(n)}) =  \left\{ \begin{array}{l} 1 \quad \mathrm{if} \quad n = k-1 \\ \sigma^k(x) \quad \mathrm{if} \quad n = k \end{array} \right. .
\end{eqnarray*}
And in the end, we get
\begin{displaymath}
\partial_{\sigma}^{[k]} \circ x = \sigma^k(x) \partial_{\sigma}^{[k]} +  \partial_{\sigma}^{[k-1]}. \quad \Box
\end{displaymath}
\end{pf}

\begin{cor} \label{invgen}
If $x \in A^\times$, then
\begin{displaymath}
\forall k \in \mathbb N \setminus \{0\}, \quad \partial_{\sigma}^{[k]} \circ x^{-1} = \sum_{i=0}^k \frac{(-1)^i}{\prod_{j=0}^i \sigma^{k-j}(x)} \partial_{\sigma}^{[k-i]}.
\end{displaymath}
\end{cor}

\begin{pf}
We compute
\begin{eqnarray*}
\sum_{i=0}^k \frac{(-1)^i}{\prod_{j=0}^i \sigma^{k-j}(x)} \partial_{\sigma}^{[k-i]} \circ x 
&=& \sum_{i=0}^k \frac{(-1)^i}{\prod_{j=0}^i \sigma^{k-j}(x)}(\sigma^{k-i}(x) \partial_{\sigma}^{[k-i]} +  \partial_{\sigma}^{[k-i-1]})
\\
&=& \sum_{i=0}^k \frac{(-1)^i}{\prod_{j=0}^{i-1} \sigma^{k-j}(x)} \partial_{\sigma}^{[k-i]} + \sum_{i=0}^{k-1} \frac{(-1)^i}{\prod_{j=0}^i \sigma^{k-j}(x)}\partial_{\sigma}^{[k-i-1]} = \partial_{\sigma}^{[k]}. \quad \Box
\end{eqnarray*}
\end{pf}

\begin{prop} \label{itform}
We have in $\mathrm D_{A/R,\sigma}^{(\infty)}$,
\begin{displaymath}
\forall k > 0, \quad \partial^{[k]}_{\sigma} \circ x = q^kx  \partial^{[k]}_{\sigma} + (k)_{q}h \partial^{[k]}_{\sigma} + \partial_{\sigma}^{[k-1]}.
\end{displaymath}
We have in $\mathrm D_{A/R,\sigma}$,
\begin{displaymath}
\forall k > 0, \quad \partial^k_{\sigma} \circ x = q^kx  \partial^{k}_{\sigma} + (k)_{q} (h \partial^{k}_{\sigma} + \partial_{\sigma}^{k-1}).
\end{displaymath}
\end{prop}
 
\begin{pf}
The first assertion is simply a reformulation of proposition \ref{commix}.
For the second one, after a base change, we may reduce to the generic case \ref{genc} and thus assume that $R$ is $q$-divisible and $q\mathrm{-char}(R) = 0$.
And we may then replace $\mathrm D_{\sigma}$ with $\mathrm D^{(\infty)}_{\sigma}$ in which case we fall back onto the first equality.
$\quad \Box$
\end{pf}

Note that, in order to prove the second formula, we cannot use directly corollary \ref{formd}: our equality takes place into the Weyl algebra and not is not just an assertion about endomorphisms of $A$.

We concentrate now on the case $\sigma(x) = qx$.

\begin{prop}
Assume $h = 0$, $q \in R^\times$ and $x \in A^\times$. 
Then, we have in $\mathrm D^{(\infty)}_{A/R,\sigma}$,
\begin{displaymath}
\forall k \in \mathbb N, \quad \partial_{\sigma}^{[k]} \circ x^{-1} = \sum_{i=0}^k (-1)^i q^{-\frac {(2k-i)(i+1)}2}x^{-i-1} \partial_{\sigma}^{[k-i]}.
\end{displaymath}
We have in $\mathrm D_{A/R,\sigma}$,
\begin{displaymath}
\forall k \in \mathbb N, \quad  \partial_{\sigma}^{k} \circ x^{-1} = \sum_{i=0}^k (-1)^i q^{-\frac {(2k-i)(i+1)}2} (k)_{q} \cdots (k-i+1)_{q} x^{-i-1} \partial_{\sigma}^{k-i}.
\end{displaymath}
\end{prop}

\begin{pf}
The first assertion is a particular case of corollary \ref{invgen} and the second one follows by the standard generic argument.
$\quad \Box$
\end{pf}

\begin{prop}
Assume that $h = 0$.
Then, we have in $\mathrm D^{(\infty)}_{A/R,\sigma}$,
\begin{displaymath}
\forall k, n \in \mathbb N, \quad \partial_{\sigma}^{[k]} \circ x^n = \sum_{i=0}^k q^{(n-i)(k-i)}{n \choose i}_{q} x^{n-i}\partial_{\sigma}^{[k-i]}.
\end{displaymath}
We have in $\mathrm D_{A/R,\sigma}$,
\begin{displaymath}
\forall k, n \in \mathbb N, \quad \partial_{\sigma}^{k} \circ x^n = \sum_{i=0}^k q^{(n-i)(k-i)}[i]_{q}!{k \choose i}_{q}{n \choose i}_{q} x^{n-i}\partial_{\sigma}^{k-i}.
\end{displaymath}
\end{prop}

\begin{pf}
As usual, the second formula will follow from the first that we prove directly by induction on $n$.
We will have
\begin{eqnarray*}
\partial_{\sigma}^{[k]} \circ x^n &=& \sum_{i=0}^k q^{(n-1-i)(k-i)}{n-1 \choose i}_{q} x^{n-1-i}\partial_{\sigma}^{[k-i]} \circ x
\\
&=& \sum_{i=0}^k q^{(n-1-i)(k-i)}{n-1 \choose i}_{q} x^{n-1-i}\left(q^{k-i}x \partial_{\sigma}^{[k-i]} + \partial_{\sigma}^{[k-i-1]}\right)
\\
&=&
\sum_{i=0}^k q^{(n-1-i)(k-i)}{n-1 \choose i}_{q} x^{n-1-i}q^{k-i} x\partial_{\sigma}^{[k-i]}
\\
&&
+ \sum_{i=0}^k q^{(n-1-i)(k-i)}{n-1 \choose i}_{q} x^{n-1-i} \partial_{\sigma}^{[k-i-1]}
\\
&=& \sum_{i=0}^k q^{(n-i)(k-i)}{n-1 \choose i}_{q} x^{n-i} \partial_{\sigma}^{[k-i]}+ \sum_{i=0}^{k-1} q^{(n-i)(k-i+1)}{n-1 \choose i-1}_{q} x^{n-i} \partial_{\sigma}^{[k-i]}
\\
&=& \sum_{i=0}^k q^{(n-i)(k-i)}\left( {n-1 \choose i}_{q} + q^{n-i}{n-1 \choose i-1}_{q}\right) x^{n-i} \partial_{\sigma}^{[k-i]}
\\
&=& \sum_{i=0}^k q^{(n-i)(k-i)}{n \choose i}_{q} x^{n-i}\partial_{\sigma}^{[k-i]}. \quad \Box
\end{eqnarray*}
\end{pf}

\begin{cor}
When $h = 0$, we have
\begin{displaymath}
\forall n \in \mathbb N, \quad \partial^{[k]}_{\sigma}(x^n) =  \left\{ \begin{array}{cl} {n \choose k}_{q} x^{n-k} & \mathrm{if} \ n \geq k \\ 0 & \mathrm{otherwise} \end{array} \right. 
\end{displaymath}
and
\begin{displaymath}
\forall n \in \mathbb N, \quad \partial^{k}_{\sigma}(x^n) =  \left\{ \begin{array}{cl} (n)_{q} \cdots (n-k+1)_{q} x^{n-k} & \mathrm{if} \ n \geq k \\ 0 & \mathrm{otherwise} \end{array} \right. .
\end{displaymath}
\end{cor}

\begin{pf}
The first assertion follows from the proposition and the second one uses the inclusion $\overline D_{A/R,\sigma} \subset D_{A/R,\sigma}^{(\infty)}$.
$\quad \Box$
\end{pf}

\begin{cor} \label{minuso}
Assume that $h = 0$.
Then, we have
\begin{displaymath}
\forall n \in \mathbb N \setminus \{0\}, \quad \partial_{\sigma} (x^{n} ) = (n)_{q} x^{n-1}.
\end{displaymath}
If $x \in A^\times$ and $q \in R^\times$, then the formula actually holds for all $n \in \mathbb Z$.
\end{cor}

\begin{pf}
Only the second part needs a proof: if $x$ is invertible, we have for $n>0$
\begin{eqnarray*}
0 &=& \partial_{\sigma}(1) = \partial_{\sigma}(x^nx^{-n}) =  \partial_{\sigma}(x^n)x^{-n} + \sigma(x^n) \partial_{\sigma}(x^{-n}) 
\\
&=& (n)_{q} x^{n-1}x^{-n} +  q^nx^n \partial_{\sigma}(x^{-n}) = (n)_{q} x^{-1} +  q^nx^n \partial_{\sigma}(x^{-n}),
\end{eqnarray*}
from which we derive, when $q \in R^\times$,
\begin{displaymath}
\partial_{\sigma}(x^{-n}) = -q^{-n}  (n)_{q} x^{-n-1} = (-n)_{q} x^{-n-1}. \quad \Box
\end{displaymath}
\end{pf}

\section{Formal deformations of twisted differential operators}

In this section, we study the relation between twisted differential operators relative to various endomorphisms of $A$.
We are particularly interested in the comparison of our twisted differential operators with usual differential operators.
We assume that there exists a twisted coordinate $x$, that we fix for the rest of the section, and we write $y_{\sigma} := x - \sigma(x)$.

Recall that the ring $\mathrm D_{\sigma}^{(\infty)}$ of differential operators (of infinite level) comes with a natural increasing filtration by $A$-submodules
\begin{displaymath}
\mathrm D_{\sigma}^{(\infty)} = \cup_{m=0}^{\infty} \mathrm{Diff}_{m,\sigma}(A),
\end{displaymath}
which is called the \emph{order filtration}.
The choice of the twisted coordinate determines a splitting of this filtration.
To see this, we let for all $m \in \mathbb N$, $K_{\sigma}^{[m]} \subset \mathrm D_{\sigma}^{(\infty)}$ be the free $A$-module generated by all $\partial_{\sigma}^{[k]}$ for $k \geq m$.
Note that this is actually a filtration by \emph{left} ideals.

\begin{dfn} \label{infinf}
The decreasing filtration by the $K_{\sigma}^{[m]}$ is called the \emph{ideal filtration} on $\mathrm D_{\sigma}^{(\infty)}$.
The \emph{module of twisted differential operators of infinite level and infinite order} on $A$ is
\begin{displaymath}
\widehat {\mathrm D}_{A/R,\sigma}^{(\infty)} = \varprojlim \mathrm D_{\sigma}^{(\infty)}/K_{\sigma}^{[m+1]}.
\end{displaymath}
\end{dfn}

We might again drop the index $A/R$ and write $\widehat {\mathrm D}_{\sigma}^{(\infty)}$.
The decreasing filtration of $\widehat {\mathrm D}_{\sigma}^{(\infty)}$ by the closures $\widehat K_{\sigma}^{[m]}$ of the $K_{\sigma}^{[m]}$'s will also be called the \emph{ideal filtration} .

\begin{rmks}
\begin{enumerate}
\item
The ideal filtration is \emph{separated}, which means that we have $\cap_{m} K_{\sigma}^{[m+1]} = \{0\}$ and it follows that $\mathrm D_{\sigma}^{(\infty)} \subset \widehat {\mathrm D}_{\sigma}^{(\infty)}$.
Actually, any $\varphi \in \widehat {\mathrm D}_{\sigma}^{(\infty)}$ can be uniquely written as an infinite sum $\sum_{0}^\infty z_{k} \partial_{\sigma}^{[k]}$ with $z_{k} \in A$ (and conversely).
In other words, we have the isomorphisms of $A$-modules
\begin{displaymath}
\mathrm D_{\sigma}^{(\infty)}= \bigoplus_{k \in \mathbb N} A \partial_{\sigma}^{[k]} \quad \mathrm{and} \quad \widehat{\mathrm D}_{\sigma}^{(\infty)} = \prod_{k \in \mathbb N} A \partial_{\sigma}^{[k]}.
\end{displaymath}
\item
The ideal filtration defines a splitting of the order filtration in the sense that
\begin{displaymath}
\mathrm D_{\sigma}^{(\infty)} = \mathrm{Diff}_{m,\sigma}(A) \oplus K_{\sigma}^{[m+1]} \quad \mathrm{and} \quad \widehat {\mathrm D}_{\sigma}^{(\infty)} = \mathrm{Diff}_{m,\sigma}(A) \oplus \widehat K_{\sigma}^{[m+1]}
\end{displaymath}
as $A$-modules.
\item
$\widehat {\mathrm D}_{\sigma}^{(\infty)}$ is \emph{not} a ring in general.
More precisely, multiplication on $\mathrm D_{\sigma}^{(\infty)}$ is not continuous for the ideal filtration:
we always have $\partial_{\sigma}^{[m]} \to 0$ when $m \to \infty$ but, if $\sigma(x) = qx$ and $x \in A^\times$ for example, we can see that
\begin{displaymath}
\forall m \in \mathbb N, \quad \partial_{\sigma}^{[m]} \circ x^{-1} \equiv \partial_{\sigma}^{[m]}(x^{-1}) = x^{-m-1} \neq 0 \mod K_{\sigma}^{[1]}.
\end{displaymath}
\end{enumerate}
\end{rmks}

If $A$ is any ring and $A[\xi]$ denotes the polynomial ring in one variable $\xi$, then the \emph{natural filtration} on $\mathrm{Hom}_{A}(A[\xi], A)$ is the decreasing filtration by the kernels of the surjections
\begin{displaymath}
\mathrm{Hom}_{A}(A[\xi], A) \to  \mathrm{Hom}_{A}(A[\xi]_{\leq m}, A),
\end{displaymath}
where $A[\xi]_{\leq m}$ denotes as before the $A$-module of polynomials of degree at most $m$.
The corresponding topology will be called the \emph{natural topology} of $\mathrm{Hom}_{A}(A[\xi], A)$.
Note that $\mathrm{Hom}_{A}(A[\xi], A)$ is separated and complete for the natural topology:
\begin{displaymath}
\varprojlim \mathrm{Hom}_{A}(A[\xi]_{\leq m}, A) = \mathrm{Hom}_{A}(\varinjlim A[\xi]_{\leq m}, A) = \mathrm{Hom}_{A}(A[\xi], A).
\end{displaymath}

\begin{lem}[Formal density] \label{density}
The map
\begin{displaymath}
\xymatrix@R0cm{ A[\xi] \ar[r] & P_{A/R}
\\ \xi \ar@{|->}[r] & \tilde x - x}
\end{displaymath}
induces by duality an isomorphism of topological $A$-modules
\begin{equation} \label{dens}
\xymatrix@R0cm{ \widehat {\mathrm D}_{A/R,\sigma}^{(\infty)} \ar[r]^-{\simeq} & \mathrm{Hom}_{A}(A[\xi], A).}
\end{equation}
More precisely, the ideal filtration corresponds to the natural filtration.
\end{lem}

\begin{pf}
It follows from corollary \ref{free} that the cokernel of the map
\begin{displaymath}
\xymatrix@R0cm{ A[\xi]_{\leq m} \ar@{^{(}->}[r] & A[\xi] \ar[r] & P \ar@{->>}[r] & P_{(n)_{\sigma}}}
\end{displaymath}
is generated by $\xi^{(k)}$ for $m < k \leq n$.
Moreover, this map is injective when $m \leq n$.
Dually, it means that the corresponding map
\begin{displaymath}
\xymatrix@R0cm{\mathrm{Diff}_{n,\sigma}(A) \simeq \mathrm{Hom}_{A}(P_{(n)_{\sigma}}, A) \ar[r] & \mathrm{Hom}_{A}(A[\xi]_{\leq m}, A)}
\end{displaymath}
is surjective for $m \leq n$ and that its kernel is exactly $K_{\sigma}^{[m+1]} \cap \mathrm{Diff}_{n,\sigma}(A)$ (it is generated by $\partial_{\sigma}^{[k]}$ for $n \geq k > m$).
Taking direct limits on the left, we see that the map
\begin{displaymath}
\xymatrix@R0cm{\mathrm D_{\sigma}^{(\infty)} \ar[r] & \mathrm{Hom}_{A}(A[\xi]_{\leq m}, A)}
\end{displaymath}
is surjective with kernel exactly $K_{\sigma}^{[m+1]}$.
In other words, we get a canonical isomorphism of $A$-modules
\begin{displaymath}
\xymatrix@R0cm{\mathrm D_{\sigma}^{(\infty)}/K_{\sigma}^{[m+1]} \ar[r]^-{\simeq} & \mathrm{Hom}_{A}(A[\xi]_{\leq m}, A)}.
\end{displaymath}
Thus, taking inverse limits on both sides gives the result.
$\quad \Box$
\end{pf}

\begin{rmks}
\begin{enumerate}
\item
By construction, there exists a commutative diagram
\begin{displaymath}
\xymatrix{\mathrm D_{A/R,\sigma}^{(\infty)} \ar@{^{(}->}[r] \ar@{^{(}->}[dd] & \mathrm{End}_{R}(A) \ar[d]^\simeq
\\
 & \mathrm{Hom}_{A}(P, A) \ar[d]
\\
\widehat {\mathrm D}_{A/R,\sigma}^{(\infty)} \ar[r]^-{\simeq} & \mathrm{Hom}_{A}(A[\xi], A)
}.
\end{displaymath}
\item
We will usually denote by $\tilde \varphi \in \mathrm{Hom}_{A}(A[\xi], A)$ the image of $\varphi \in \widehat {\mathrm D}_{A/R,\sigma}^{(\infty)}$.
This is compatible with the previous notation for linearization.
In particular, we have
\begin{displaymath}
\forall k, n \in \mathbb N, \quad \widetilde{\partial^{[k]}_{\sigma}}(\xi^{(n)}) = \left\{\begin{array} {cll} 1 & \mathrm{if} \ k = n \\ 0 & \mathrm{otherwise}\end{array} \right. .
\end{displaymath}
\item
When $A = R[x]$, we actually get
\begin{displaymath}
\xymatrix@R0cm{\mathrm D_{A/R,\sigma}^{(\infty)} \ar@{^{(}->}[r] & \widehat {\mathrm D}_{A/R,\sigma}^{(\infty)} \ar[r]^-{\simeq} & \mathrm{Hom}_{A}(A[\xi], A) \ar[r]^-{\simeq} & \mathrm{Hom}_{A}(P, A) \ar[r]^-{\simeq} & \mathrm{End}_{R}(A)}
\end{displaymath}
and $\widehat {\mathrm D}_{A/R,\sigma}^{(\infty)}$ is a ring in this very special case.
\item There exists a complex analytic analog of the density lemma (for usual differential operators) as explained by Zogman Mebkhout and Luis Narv\'aez in \cite{MebkhoutNarvaez98}:
the ring of algebraic differential operators is dense in the ring of continuous endomorphisms of the structural sheaf.
\end{enumerate}
\end{rmks}

Although it is difficult to give an explicit description of the isomorphism \eqref{dens}, we can at least show the following:

\begin{lem} \label{compone}
We have 
\begin{displaymath}
\forall n \in \mathbb N \setminus \{0\}, \quad \tilde \partial_{\sigma}(\xi^n) = (\sigma(x) - x)^{n-1}.
\end{displaymath}
More generally, if $\tau$ is another $R$-algebra endomorphism of $A$, we have
\begin{displaymath}
\forall n \in \mathbb N \setminus \{0\}, \quad \tilde \partial_{\sigma}(\xi^{(n)_{\tau}}) = \prod_{i=1}^{n-1}(\sigma(x) - \tau^i(x)),
\end{displaymath}
where $\xi^{(n)_{\tau}} := \xi\tau(\xi) \cdots \tau^{n-1}(\xi)$ denotes the twisted power with respect to the endomorphism $\tau$.
\end{lem}

\begin{pf}
By definition, $\tilde \partial_{\sigma}$ is the unique $A$-linear function on $A[\xi]$ such that $\tilde \partial_{\sigma}(\xi^{(n)}) = 1$ for $n = 1$ and $0$ otherwise.
We consider now the unique $A$-linear function $u$ on $A[\xi]$ such that $u(\xi^n) = (\sigma(x) - x)^{n-1}$ for $n > 0$ and $u(1) = 0$.
We want to show that $u = \tilde \partial_{\sigma}$.
The map $u$ may be seen as the composition of division by $\xi$ on $A[\xi]$ (after removing the constant term) and evaluation at $\sigma(x) -x$.
But we have
\begin{displaymath}
\xi^{(n)} := \xi\sigma(\xi) \cdots \sigma^{n-1}(\xi),
\end{displaymath}
and we know that $\sigma(\xi) = \xi +x - \sigma(x)$.
Therefore, it is clear that for $n \geq 2$, we will get $u(\xi^{(n)}) = 0$.
Of course, for $n = 0$ we have $\xi^{(n)} = 1$ and we also obtain $0$.
Finally, for $n = 1$, we have $\xi^{(n)} = \xi$ and we get $1$.

The proof of the second formula follows the same lines.
From the first part, we can interpret $\tilde \partial_{\sigma}$ as the composition of division by $\xi$ on $A[\xi]$ and evaluation at $\sigma(x) -x$.
We want to apply this to $\xi^{(n)_{\tau}} = \prod_{i=0}^{n-1} \tau^i(\xi)$ and we know that $ \tau^i(\xi) = \xi + x - \tau^i(x)$.
Thus, we see that
\begin{displaymath}
\tilde \partial_{\sigma}(\xi^{(n)_{\tau}}) =  \prod_{i=1}^{n-1} (\sigma(x) -x + x - \tau^i(x)) = \prod_{i=1}^{n-1}(\sigma(x) - \tau^i(x)). \quad \Box
\end{displaymath}
\end{pf}

We can actually derive a remarkable consequence of the density lemma:

\begin{prop}[Formal deformation] \label{dfthm}
If $x$ is a also a $\tau$-coordinate for another $R$-endomorphism $\tau$ of $A$, then there exists an isomorphism of topological $A$-modules that depends only on $x$
\begin{displaymath}
\xymatrix@R0cm{ \widehat {\mathrm D}_{A/R,\sigma}^{(\infty)} \ar[r]^-{\simeq} & \widehat {\mathrm D}_{A/R,\tau}^{(\infty)}}. 
\end{displaymath}
More precisely, it is compatible with the ideal filtrations.
\end{prop}

Recall that the hypothesis is always satisfied when $A$ is a $\tau$-twisted localization of $R[x]$.

\begin{pf}
We may just compose the isomorphism of the density lemma with the inverse of the analogous isomorphism for $\tau$.
$\quad \Box$
\end{pf}

When the coordinate $x$ is fixed, we may (and will) \emph{identify} these two topological (or even, filtered) $A$-modules.
It is worth mentioning a particular case of this corollary that is of great interest (this is the case $\tau = \mathrm{Id}_{A})$:

\begin{cor} \label{compid}
Assume that $x$ is also an usual coordinate on $A$.
Then, there exists a canonical isomorphism
\begin{displaymath}
\xymatrix@R0cm{ \widehat {\mathrm D}_{A/R,\sigma}^{(\infty)} \ar[r]^-{\simeq} & \widehat {\mathrm D}_{A/R}^{(\infty)}.} \quad \Box
\end{displaymath}
\end{cor}

\begin{rmk}
In  theorem 2 of the introduction to \cite{Pulita14}, Pulita shows that, in the $p$-adic world, some differential modules have a natural structure of $\sigma$-module.
The main idea is to realize a formal solution of the differential module as formal solution of some $\sigma$-module.
This is analog to the way we derive the formal deformation theorem from the density lemma.
It seems that the first result in this direction was obtained by Andr\'e and Di Vizio in \cite{AndreDiVizio04}.
\end{rmk}

We can give explicit formulas in order to express the twisted derivation from one world as a twisted differential operators in another world:

\begin{prop} \label{changsig}
Assume that $x$ is a also a $\tau$-coordinate for some other $R$-endomorphism $\tau$ of $A$.
Then, we have
 \begin{displaymath}
\partial_{\sigma} = \sum_{k=1}^{\infty} \left(\prod_{i=1}^{k-1}(\sigma(x) - \tau^i(x))\right)\partial_{\tau}^{[k]}.
\end{displaymath}
\end{prop}

\begin{pf}
This follows directly by duality from lemma \ref{compone}.
$\quad \Box$
\end{pf}

\begin{cor} \label{siginf} \label{corsig}
Assume that $x$ is also a usual coordinate on $A$.
Then,
\begin{equation*} \label{transform}
\partial_{\sigma} = \sum_{k=1}^{\infty} (\sigma(x) -x)^{k-1} \partial^{[k]}
\end{equation*}
and
\begin{equation*} \label{transform2}
\sigma = \sum_{k=0}^{\infty} (\sigma(x) -x)^{k} \partial^{[k]}.
\end{equation*}
\end{cor}

\begin{pf}
The first assertion follows directly from the proposition and the other one is deduced from the equality $\sigma = 1 - y_{\sigma}\partial_{\sigma}$.
$\quad \Box$
\end{pf}

\begin{rmks}
\begin{enumerate}
\item
In the case $\sigma(x) = qx$, the formulas read
\begin{displaymath}
\partial_{\sigma} = \sum_{k=1}^{\infty} (q-1)^{k-1}x^{k-1} \partial^{[k]} \quad \mathrm{and} \quad \sigma = \sum_{k=0}^{\infty} (q-1)^{k}x^{k} \partial^{[k]}.
\end{displaymath}
\item
In the case $\sigma(x) = x + h$, the formulas read
\begin{displaymath}
\partial_{\sigma} = \sum_{k=1}^{\infty} h^{k-1} \partial^{[k]} \quad \mathrm{and} \quad \sigma = \sum_{k=0}^{\infty} h^{k} \partial^{[k]}.
\end{displaymath}
\end{enumerate}
\end{rmks}

Conversely, we may also express the usual derivation in term of twisted derivations.
We only do the classical cases:

\begin{cor}
Assume that $x$ is also a usual coordinate on $A$.
Then,
\begin{enumerate}
\item
If $\sigma(x) = qx$, we have
\begin{displaymath}
\partial = \sum_{k=1}^{\infty} (1-q)^{k-1}(k-1)_{q}! x^{k-1}\partial_{\sigma}^{[k]}.
\end{displaymath}
\item
If $\sigma(x) = x + h$, we have
\begin{displaymath}
\partial = \sum_{k=1}^{\infty} (-1)^{k-1}h^{k-1}(k-1)!\partial_{\sigma}^{[k]}.
\end{displaymath}

\end{enumerate}
\end{cor}

\begin{pf}
 In the formula of proposition \ref{changsig} we replace $\tau$ by $\sigma$ and $\sigma$ by $\mathrm{Id}$ respectively, in order to obtain
\begin{eqnarray*}
\partial &=& \sum_{k=1}^{\infty} \prod_{i=1}^{k-1}(x -q^ix)\partial_{\sigma}^{[k]} = \sum_{k=1}^{\infty}\prod_{i=1}^{k-1}(1 -q^i)x^{k-1}\partial_{\sigma}^{[k]}
\\
&=& \sum_{k=1}^{\infty}(1-q)^{k-1}\prod_{i=1}^{k-1}(i)_{q}x^{k-1}\partial_{\sigma}^{[k]} = \sum_{k=1}^{\infty} (1-q)^{k-1}(k-1)_{q}! x^{k-1}\partial_{\sigma}^{[k]}
\end{eqnarray*}
in the first case, and
\begin{displaymath}
\partial = \sum_{k=1}^{\infty} \prod_{i=1}^{k-1}(-ih)\partial_{\sigma}^{[k]} = \sum_{k=1}^{\infty} (-h)^{k-1}(k-1)!\partial_{\sigma}^{[k]}
\end{displaymath}
in the other one.
$\quad \Box$
\end{pf}

We may also apply proposition \ref{changsig} in order to make more precise the statement of proposition \ref{pow}.

\begin{cor} \label{prodf}
Assume that $x$ is also a $\sigma^m$-coordinate for $A$ over $R$, then we have
\begin{displaymath}
\partial_{\sigma^m} = \sum_{k=1}^{m} \left(\prod_{i=1}^{k-1}(\sigma^m(x) - \sigma^i(x))\right)\partial_{\sigma}^{[k]}. \quad \Box
\end{displaymath}
\end{cor}

One can give a more concrete formula in the quantum situation:

\begin{cor} \label{bigform}
Assume that $x$ is also a $\sigma^m$-coordinate for $A$ over $R$ and that $\sigma(x) = qx$.
Then, we have
\begin{displaymath}
\partial_{\sigma^m} = \sum_{k=1}^{m} q^{\frac {k(k-1)}2}(q - 1)^{k-1} (m-1)_{q}\cdots (m-k+1)_{q} x^{k-1} \partial_{\sigma}^{[k]}.
\end{displaymath}
\end{cor}

\begin{pf}
We compute the coefficient of $\partial_{\sigma}^{[k]}$ for $1 \leq k \leq m$:
\begin{eqnarray*}
\prod_{i=1}^{k-1}(q^mx - q^ix)
&=&
\left(\prod_{i=1}^{k-1} q^i\right)\left(\prod_{i=1}^{k-1}(q^{m-i} - 1)\right)x^{k-1}
\\
&=& q^{\frac {k(k-1)}2}(q - 1)^{k-1}\left(\prod_{i=1}^{k-1} (m-i)_{q}\right)x^{k-1}. \quad \Box
\end{eqnarray*}
\end{pf}

\begin{xmp}
We will have $
\partial_{\sigma^2} = \partial_{\sigma} + q(q-1)x\partial^{[2]}_{\sigma}.
$
\end{xmp}

\section{Formal confluence}

We explain here how one can use the quantum Weyl algebra in order to approximate a usual differential operator.
We assume that $A$ is a quantum $R$-algebra which means that we are given a twisted coordinate $x$ on $A$ with $\sigma(x) = qx + h$ with $q,h \in R$.

Recall from \cite{LeStumQuiros15a} that the twisted Weyl $R$-algebra $\mathrm D_{\sigma}$ has an increasing filtration by free $A$-modules of finite rank
\begin{displaymath}
\mathrm{Fil}^m\mathrm D_{\sigma} = \bigoplus_{k\leq m} A \partial_{\sigma}^k
\end{displaymath}
that is called the \emph{order filtration}.
But it also has a decreasing filtration by free $A$-modules (of infinite rank)
\begin{displaymath}
K_{\sigma}^{m} = \bigoplus_{k\geq m} A \partial_{\sigma}^k
\end{displaymath}
that we called the \emph{ideal filtration}.
We will consider the completion
\begin{displaymath}
\widehat {\mathrm D}_{\sigma} = \varprojlim \mathrm D_{\sigma}/K_{\sigma}^{m+1}
\end{displaymath}
that also comes with its \emph{ideal filtration} and we have $\mathrm D_{\sigma} \subset \widehat {\mathrm D}_{\sigma}$.
Note that $\widehat {\mathrm D}_{\sigma}$ is the set of all infinite sums $\sum_{0}^\infty z_{k} \partial_{\sigma}^{k}$ with $z_{k} \in A$.
In other words, there exists an isomorphism of $A$-modules
\begin{displaymath}
\widehat{\mathrm D}_{\sigma} = \prod_{k \in \mathbb N} A \partial_{\sigma}^k.
\end{displaymath}

The $A$-module $\widehat {\mathrm D}_{\sigma}$ is not a ring in general.
However, the result holds in the finite quantum characteristic case as we can check right now:
 
\begin{prop}
If $A$ is $q$-flat and $q\mathrm{-char}(A) = p > 0$, then multiplication is continuous for the ideal topology on $\mathrm D_{A/R,\sigma}$ and turns $\widehat{\mathrm D}_{A/R,\sigma}$ into an $R$-algebra.
\end{prop}

\begin{pf}
From the equality (see the remark following theorem \ref{thd})
\begin{displaymath}
\mathrm D_{\sigma}/\partial_{\sigma}^p \simeq \overline{\mathrm D}_{\sigma},
\end{displaymath}
and the fact that $\overline{\mathrm D}_{\sigma}$ is a ring, we deduce that the (left) ideal generated by $\partial_{\sigma}^p$ is a two sided ideal.
Since multiplication is automatically continuous for the $\partial_{\sigma}^p$-adic filtration, it follows that it is also continuous for the ideal filtration (which is the filtration by left ideals generated by the powers of $\partial_{\sigma}$).
$\quad \Box$
\end{pf}

\begin{lem} \label{gdca}
The composite map
\begin{displaymath}
\xymatrix@R0cm{ \mathrm \mathrm D_{A/R,\sigma} \ar@{->>}[r] & \overline{\mathrm D}_{A/R,\sigma} \ar@{^{(}->}[r] & \mathrm D^{(\infty)}_{A/R,\sigma}
}
\end{displaymath}
is compatible with the ideal filtrations.
Moreover, if $R$ is $q$-divisible and $q\mathrm{-char}(R) = 0$, then $\widehat {\mathrm D}_{A/R,\sigma} \simeq \widehat {\mathrm D}_{A/R,\sigma}^{(\infty)}$.
\end{lem}

This applies in particular when $R$ is a $\mathbb Q$-algebra and $\sigma = \mathrm{Id}$.

\begin{pf}
The first assertion follows from corollary \ref{formd} and the second from theorem \ref{thd}.
$\quad \Box$
\end{pf}

We can now state our fist confluence theorem:

\begin{thm}[Formal quantum confluence 1] \label{confalg}
Let $R$ be a $\mathbb Q$-algebra, $A$ a $q$-$R$-algebra for some $q \in R$.
Assume that the quantum coordinate on $A$ is also a usual coordinate.
Then, there exists a canonical map of filtered $A$-modules
\begin{equation} \label{confone}
\mathrm D_{A/R,\sigma} \to \widehat {\mathrm D}_{A/R}.
\end{equation}
Moreover, if $R$ is $q$-divisible, we have
\begin{enumerate}
\item If $q\mathrm{-char}(A) = 0$, then the map \eqref{confone} is injective and $\mathrm D_{A/R,\sigma}$ is dense in $\widehat {\mathrm D}_{A/R}$.
Actually, the inclusion $\mathrm D_{A/R,\sigma} \subset \widehat {\mathrm D}_{A/R}$ is strict for the ideal filtrations.
\item If $q\mathrm{-char}(A) = p > 0$, then the map \eqref{confone} induces an isomorphism
\begin{displaymath}
\left(\mathrm D_{A/R,\sigma}/\partial_{\sigma}^p \simeq \right) \quad \overline {\mathrm D}_{A/R,\sigma} \simeq \mathrm D_{A/R}/\partial^p.
\end{displaymath}
\end{enumerate}
\end{thm}

\begin{pf}
The map \eqref{confone} is simply the composite
\begin{displaymath}
\mathrm D_{\sigma} \to \mathrm D_{\sigma}^{(\infty)} \hookrightarrow \widehat{\mathrm D}_{\sigma}^{(\infty)} \simeq \widehat{\mathrm D}^{(\infty)} \simeq \widehat {\mathrm D}_{A/R},
\end{displaymath}
where the next to the last map is the formal deformation isomorphism of corollary \ref{compid} and the last one comes from lemma \ref{gdca} applied to the case $\sigma = \mathrm{id}$ since $R$ is a $\mathbb Q$-algebra.

If we assume that $R$ is $q$-divisible and that $q\mathrm{-char}(A) = 0$, then lemma \ref{gdca} tells us that $\mathrm D_{\sigma} = \mathrm D^{(\infty)}_{\sigma}$ as filtered rings (for the ideal filtrations) and we can use corollary \ref{compid} again.

Finally, when $R$ is $q$-divisible but $q\mathrm{-char}(A) = p >0$, then $(p)_{q} = 0$ in $A$ and all $(m)_{q} \in R^\times$ for $m < p$.
We can use the last assertion of theorem \ref{thd} and the fact that the isomorphism of corollary \ref{compid} is strictly compatible with the filtrations.
$\quad \Box$
\end{pf}

As in \cite{LeStumQuiros15a}, we denote by $A[T]_{\sigma}$ the twisted polynomial ring associated to $\sigma$: this is the non commutative polynomial ring in $T$ over $A$ with the commutation rule $Tx = \sigma(x)T$.
Recall also that the twisted coordinate $x$ is said to be \emph{strong} is $x - \sigma(x) \in A^\times$.

\begin{cor}
If $R$ is a $q$-divisible $\mathbb Q$-algebra, $q\mathrm{-char}(A) = 0$, $x$ is strong and is also a coordinate for all $\sigma^n$'s, then, the $A$-linear map
\begin{equation} \label{confilm}
\xymatrix@R0cm{A[T]_{\sigma} \ar[r] & \widehat{\mathrm D}_{A/R}
\\ T^{n} \ar@{|->}[r] & \sum_{k=0}^{\infty} \frac 1{k!}(\sigma^n(x) -x)^{k} \partial^{k}}
\end{equation}
has dense image.
\end{cor}

\begin{pf}
Compose the isomorphism of theorem 6.13 of \cite{LeStumQuiros15a} with the map of theorem \ref{confalg}.
The formula comes from corollary \ref{siginf}.
$\quad \Box$
\end{pf}

\begin{xmp}
The theorem (and its corollary) applies in particular when $q$ is not a root of unity in some field $K$ of characteristic zero, $R$ is an algebra containing $K$, $A = R[x, x^{-1}]$ and $\sigma(x) = qx$.
In this situation, $\mathrm D := \mathrm D_{A/R}$ is the non commutative ring $R[x, x^{-1}, \partial]$ with the commutation rule $\partial x = x \partial +1$ and we may see $\widehat{\mathrm D}$ (which is \emph{not} a ring) as the set of infinite sums $R[x, x^{-1}][[\partial]]$.
Also, $\mathrm D_{\sigma} := \mathrm D_{A/R,\sigma}$ is the non commutative ring $R[x, x^{-1}, \partial_{\sigma}]$ with the commutation rule $\partial_{\sigma} x = qx \partial_{\sigma} +1$.
The map of the theorem satisfies
\begin{displaymath}
\partial_{\sigma} \mapsto \sum_{k=1}^{\infty} \frac 1 {k!}(q-1)^{k-1}x^{k-1} \partial^{k},
\end{displaymath}
and the map of the corollary is given by
\begin{displaymath}
T^{n} \mapsto \sum_{k=0}^{\infty} \frac 1 {k!}(q^{n} -1)^{k} x^{k}\partial^{k}.
\end{displaymath}
\end{xmp}

In the next section, we will have to move between $\sigma$ and powers of $\sigma$ (or more precisely, the other way around, between $\sigma$ and roots of $\sigma$).

\begin{prop} \label{explem}
Let $m \in \mathbb N \setminus \{0\}$.
If $(m)_{q}! \in R^\times$, then there exists a unique $A$-linear ring homomorphism
\begin{displaymath}
\xymatrix@R0cm{\mathrm D_{A/R,\sigma^m} \ar[r]^{\iota_{\sigma, m}} & \mathrm D_{A/R,\sigma}
\\ \partial_{\sigma^m} \ar@{|->}[r] & \sum_{k=1}^{m} \frac 1{(k)_{q}!}\left(\prod_{i=1}^{k-1}(\sigma^m(x) - \sigma^i(x))\right)\partial_{\sigma}^{k}
}.
\end{displaymath}
Moreover, the diagram
\begin{equation} \label{iotdiag}
\xymatrix{\mathrm D_{A/R,\sigma^m} \ar[r]\ar[d]^{\iota_{\sigma, m}} & \mathrm D_{A/R,\sigma^m}^{(\infty)} \ar@{^{(}->}[d] \ar@{^{(}->}[r] & \mathrm {End}_{R}(A) \ar@{=}[d]
\\ \mathrm D_{A/R,\sigma} \ar[r]  & \mathrm D_{A/R,\sigma}^{(\infty)} \ar@{^{(}->}[r] & \mathrm {End}_{R}(A)
}
\end{equation}
is commutative.
\end{prop}

\begin{pf}
We may assume (see the second remark following theorem \ref{thd}) that we are in the generic situation
\begin{displaymath}
R = \mathbb Q(t)[s], \quad A = R[x], \quad q = t \quad \mathrm{and} \quad h = s.
\end{displaymath}
But then, the left horizontal arrows in diagram \eqref{iotdiag} are bijective.
In other words, we can identify the twisted Weyl $R$-algebras with the rings of twisted differential operators of infinite level and use corollary \ref{prodf} and corollary \ref{formd}.
\end{pf}

\begin{rmks}
\begin{enumerate}
\item
In the case $h = 0$, we have thanks to corollary \ref{bigform}, the more concrete formula:
\begin{displaymath}
\xymatrix@R0cm{\partial_{\sigma^m} \ar@{|->}[r] & \sum_{k=1}^{m} \frac { q^{\frac {k(k-1)}{2}}(q - 1)^{k-1}} {(k)_{q}} {m-1 \choose k-1}_{q}x^{k-1}  \partial_{\sigma}^{k}.
}
\end{displaymath}
\item
Assume that $R$ is $q$-divisible.
If $q\mathrm{-char}(R) = p > 0$, then the hypothesis is satisfied if and only if $m < p$.
Of course if $q\mathrm{-char}(R) = 0$, then the hypothesis is satisfied for any $m \in \mathbb N$.
\item
\end{enumerate}
\end{rmks}

\begin{prop}
Let $m,n \in \mathbb N \setminus \{0\}$.
If $(mn)_{q}!$ is invertible, then
\begin{displaymath}
\iota_{\sigma, mn} = \iota_{\sigma^m,n} \circ \iota_{\sigma,m}.
\end{displaymath}
\end{prop}

\begin{pf}
After a base change, we may assume that all $q$-integers are invertible and identify twisted Weyl $R$-algebras with rings of twisted differential operators of infinite level.
The assertion then becomes a triviality.
$\quad \Box$
\end{pf}

\begin{prop} \label{compuiss}
If $(m)_{q}!$ is invertible in $R$, then the diagram
\begin{displaymath}
\xymatrix{T \ar@{|->}[d]  & A[T]_{\sigma^m} \ar[r] \ar[d] & \mathrm D_{A/R,\sigma^m} \ar[d]^{\iota_{\sigma, m}}, & T \mapsto 1 + (\sigma^m(x) -x)\partial_{\sigma^m}
\\ \mathrm T^m & A[T]_{\sigma} \ar[r] & \mathrm D_{A/R,\sigma} , & T \mapsto 1 + (\sigma(x) -x)\partial_{\sigma}
}
\end{displaymath}
is commutative.
\end{prop}

\begin{pf}
Again, we may assume that all $q$-integers are invertible in $R$ and use the second remark following theorem 6.13 of \cite{LeStumQuiros15a}.
$\quad \Box$
\end{pf}

\section{Formal confluence in positive quantum characteristic}

In this section, we extend the formal confluence theorem to the case of positive quantum characteristic.
In order to do that, it is necessary to use the $S$-twisted theory where $S$ a filtering (for division) set of positive integers.
Thus, we assume here that $A$ is an $S$-twisted $R$-algebra in the sense of \cite{LeStumQuiros15a} (it is endowed with a compatible system of $n$-th roots $\sigma_{n}$ of $\sigma$ for all $n \in S$).

We recall from \cite{LeStumQuiros15} that a \emph{system of roots} of $q \in R$ is a compatible family $\underline q := \{q_{n}\}_{n \in S}$ of $n$-th roots of $q$.
We call the system $\underline q$ \emph{admissible} if
\begin{displaymath}
\forall n \in S, \quad (n)_{q^n} \in R^\times.
\end{displaymath}
This is a natural condition in order to define the \emph{$q$-rational number $(r)_{q}$ for} any $r \in \mathbb N\frac 1S$.
More precisely, if $r = \frac mn$ with $m \in \mathbb N$ and $n \in S$, then,
\begin{displaymath}
(r)_{q} := \frac {(m)_{q_{n}}}{(n)_{q_{n}}}
\end{displaymath}
only depends on $r$ and not on the choice of $m$ and $n$.
It is convenient to introduce the following terminology:

\begin{dfn} Let $\{q_{n}\}_{n \in S}$ be a system of roots of $q$ in $R$.
Then, $R$ is said to be \emph{$\underline q$-divisible} if $R$ is $q_{n}$-divisible for all $n \in S$.
\end{dfn}

When the system is admissible, there exists a nice equivalent definition:

\begin{lem} Let $\{q_{n}\}_{n \in S}$ be an admissible system of roots of $q$.
Then $R$ is $\underline q$-divisible if and only if
\begin{displaymath}
\forall r \in \mathbb N\frac 1S, \quad (r)_{q} \in R^\times \cup \{0\}.
\end{displaymath}
\end{lem}

\begin{pf}
If $r = \frac mn$ with $m \in \mathbb N$ and $n \in S$, we have $(r)_{q} = \frac {(m)_{q_{n}}}{(n)_{q_{n}}}$.
It follows that $(r)_{q} \in R^\times$ (resp. $= 0$) if and only if $(m)_{q_{n}} \in R^\times$ (resp. $= 0$).
$\quad \Box$
\end{pf}

\begin{rmks}
\begin{enumerate}
\item If $1 - q \in R^\times$, then $\underline q$ is admissible.
In particular, if $q \in K$ where $K$ is a subfield of $R$ and $q \neq 0$, then $\underline q$ is admissible.
\item If $\underline q \subset K$ where $K$ is a subfield of $R$, then $R$ is $\underline q$-divisible.
\end{enumerate}
\end{rmks}

\begin{lem} For a system $\underline q := \{q_{n}\}_{n \in S}$ of roots $q$, the following are equivalent:
\begin{enumerate}
\item $\underline q$ is admissible and for all $r \in \mathbb N\frac 1S \cap [0, 1]$, we have $(r)_{q} \in R^\times$.
\item  For all $n \in S$, we have $(n)_{q_{n}}! \in R^\times$.
\end{enumerate}
\end{lem}

\begin{pf}
First of all, the second condition also implies that $\underline q$ is admissible and we may therefore make this assumption.
If $r = \frac mn$ with $m \in \mathbb N$ and $n \in S$, we know that $(r)_{q} \in R^\times$ if and only if $(m)_{q_{n}} \in R^\times$.
Also, we have $r \leq 1$ if and only $m \leq n$.
Thus, the second condition, which means that $(m)_{q_{n}} \in R^\times$ whenever $m \leq n$, is equivalent to the requirement that $(r)_{q} \in R^\times$ for $r \leq 1$.
$\quad \Box$
\end{pf}

\begin{dfn}
Such a system of roots will be called \emph{strongly admissible}.
\end{dfn}

\begin{lem} \label{charcri}
If $\underline q$ is an admissible (resp. a strongly admissible) system of roots, and we write for all $n \in S$, $p_{n} := q_{n}\mathrm{-char}(R)$, then we have
\begin{equation*} \label{carad}
\forall n \in S, \quad p_{n} \nmid n \ (\mathrm{resp.}\ p_{n} > n) \quad \mathrm{or} \quad  p_{n} = 0.
\end{equation*}
When $R$ is $\underline q$-divisible, the converse is true.
\end{lem}

\begin{pf}
By definition, $q_{n}$ is admissible (resp. strongly admissible) if and only if $(n)_{q_{n}} \in R^\times$ (resp. $(n)_{q_{n}}! \in R^\times$).
By definition also, we have $p_{n} \neq 0$ and $p_{n} \mid n$ (resp. $p_{n} \leq n$) if and only if $(n)_{q_{n}} = 0$ (resp. $(n)_{q_{n}}! = 0$).
We see that these two conditions are mutually exclusive in general and that they are exhaustive when $R$ is $q_{n}$-divisible.
$\quad \Box$
\end{pf}

There usually exists strongly admissible systems as the next result shows:

\begin{prop}
Assume that $R$ is $\underline q$-divisible and that $q\mathrm{-char}(R) = p > 0$.
Then, if $\underline q$ is a system of $p^n$-th root of $q$, it is strongly admissible.
\end{prop}

\begin{pf}
It follows from proposition 1.16 of \cite{LeStumQuiros15} that for all $n \in S$, we have $q_{p^n}\mathrm{-char}(R) = p^{n+1} > p^n$ and we can apply lemma \ref{charcri}.
$\quad \Box$
\end{pf}

\begin{dfn}
An $x \in A$ is called an \emph{$S$-twisted coordinate}, or \emph{rooted twisted coordinate}, (resp. an \emph{$S$-quantum coordinate}, or \emph{rooted quantum coordinate},) if $x$ is a twisted (resp. quantum) coordinate for all $(A, \sigma_{n})$.
We call it \emph{strong} if $x - \sigma_{n}(x) \in A^\times$ for all $n \in S$.
\end{dfn}

Thus, by definition, $x$ is an $S$-quantum coordinate if and only if
\begin{equation} \label{sigen}
\forall n \in S, \quad \sigma_{n}(x) = q_{n} x + h_{n} \quad \mathrm{with} \quad q_{n}, h_{n} \in R.
\end{equation}
When this is the case, we might also say \emph{$\underline q$-coordinate} and call $A$ a \emph{rooted quantum $R$-algebra}, an \emph{$S$-quantum $R$-algebra} or a \emph{$\underline q$-$R$-algebra.}
We call it \emph{strong} when $x$ is strong.

\begin{lem}
Assume that $x$ is simultaneously a twisted coordinate and a rooted quantum coordinate on $A$ so that \eqref{sigen} holds.
Then, $x$ is a quantum coordinate on $A$ and if we write $\sigma(x) = qx + h$ with $q, h \in R$, then $\underline q := \{q_{n}\}_{n \in S}$ is a system of roots of $q$.
If this system is admissible, we have
\begin{displaymath}
\forall n \in \mathbb N, \quad h_{n} = \left(\frac 1n\right)_{q}h.
\end{displaymath}
\end{lem}

\begin{pf}
If we are given $m \in \mathbb N$ and $n \in S$, we know that $\sigma_{n}^m$ only depends on $r := \frac mn$ and so does $\sigma_{n}^m(x) = q_{n}^mx + (m)_{q_{n}}h_{n}$.
It immediately follows that $\underline q$ is a system of roots.
Moreover, the case $m = n$ implies that  $\sigma(x) = qx + h$ with $q = q_{n}^n$ and $h = (n)_{q^n}h_{n}$ from which we derive the other assertions.
$\quad \Box$
\end{pf}

\begin{xmps}
\begin{enumerate}
\item If we are given a system $\underline q$ of roots in $R$, we can endow $A := R[x]$ or $A = R[x, x^{-1}]$ with
\begin{displaymath}
\forall n \in S, \quad \sigma_{n}(x) = q_{n} x.
\end{displaymath}
In the later case, $x$ is a strong $\underline q$-coordinate if and only if $1 -q \in R^\times$ (and then $\underline q$ is also admissible).
\item More generally, if we are given an \emph{admissible} system $\underline q$ of roots of $q \in R$ and some $h \in R$, we can endow $R[x]$ with
\begin{displaymath}
\forall n \in S, \quad \sigma_{n}(x) = q_{n} x +\left(\frac 1n\right)_{q}h.
\end{displaymath}
\end{enumerate}
\end{xmps}

\begin{dfn}
Assume that $x$ is a $\underline q$-coordinate on $A$ with $\underline q$ strongly admissible.
A \emph{rooted twisted differential $A$-module} is an $A$-module $M$ endowed with a family of $\partial_{A,\sigma_{n}}$-derivations $\partial_{M,n}$ for all $n \in S$ such that whenever $n,n' \in S$ with $n | n'$, we have
\begin{displaymath}
\forall s \in M, \quad \partial_{M,n} (s) =  \sum_{k=1}^{n'/n} \frac 1{(k)_{q_{n'}}!}\left(\prod_{i=1}^{k-1}(\sigma_{n}(x) - \sigma_{n'}^i(x))\right)\partial_{M,n'}^{k}(s).
\end{displaymath}
\end{dfn}

They form a category that we will denote by $\partial_{A,\underline \sigma}\mathrm{-Mod}^{\mathrm{root}}$.

\begin{dfn}
If $A$ is a $\underline q$-$R$-algebra with $\underline q$ strongly admissible,
then the \emph{rooted twisted Weyl algebra} of $A$ is
\begin{displaymath}
\mathrm D_{A/R,\underline \sigma}:= \varinjlim \mathrm D_{A/R,\sigma_{n}}
\end{displaymath}
with transition maps (coming from lemma \ref{explem}) for $n,n' \in S$ with $n | n'$ given by
\begin{displaymath}
\xymatrix@R0cm{\mathrm D_{A/R,\sigma_{n}} \ar[r]^{\iota_{\sigma_{n'}, n'/n}} & \mathrm D_{A/R,\sigma_{n'}}
\\ \partial_{\sigma_{n}} \ar@{|->}[r] & \sum_{k=1}^{n'/n} \frac 1{(k)_{q_{n'}}!}\left(\prod_{i=1}^{k-1}(\sigma_{n}(x) - \sigma_{n'}^i(x))\right)\partial_{\sigma_{n'}}^{k}
}.
\end{displaymath}

\end{dfn}

\begin{prop} \label{abag}
Assume that $A$ is a $\underline q$-$R$-algebra with $\underline q$ strongly admissible.
Then the rooted twisted differential $A$-modules form an abelian category with sufficiently many injective and projective objects.
Actually, if $M$ is a $\mathrm D_{A/R,\underline \sigma}$-module, then the maps 
\begin{displaymath}
\xymatrix@R0cm{ M  \ar[rr]^{\partial_{M,n}} && M 
\\ s \ar@{|->}[rr] && \partial_{M,n} s}
\end{displaymath}
turn $M$ into a rooted twisted differential $A$-module and we obtain an equivalence (an isomorphism) of categories
\begin{displaymath}
\mathrm D_{A/R,\underline \sigma}\mathrm{-Mod} \simeq \partial_{A,\underline \sigma}\mathrm{-Mod}^{\mathrm{root}}.
\end{displaymath}
\end{prop}

\begin{pf} This is a consequence of proposition 5.6 of \cite{LeStumQuiros15a}.
$\quad \Box$
\end{pf}

Recall from \cite{LeStumQuiros15a} that $A[T^{\frac 1S}]_{\underline \sigma}$ denotes the non commutative Puiseux polynomial ring with the commutation rule
\begin{displaymath}
\forall m \in \mathbb N, \forall n \in S, \quad T^{\frac mn}x = \sigma_{n}^m(x) T^{\frac mn}.
\end{displaymath}
One can also consider the notion of \emph{$\underline \sigma_{A}$-module}: this is an $A$-module endowed with a compatible family of $\sigma_{A,n}$-linear endomorphisms $\sigma_{M,n}$.
There exists an equivalence (an isomorphism) between the category of $A[T^{\frac 1S}]_{\underline \sigma}$-modules and the category of $\underline \sigma_{A}$-modules.

\begin{prop}
Assume that $x$ is a $\underline q$-coordinate on $A$ with $\underline q$ strongly admissible.
Then there exists a unique $A$-linear homomorphism of rings
\begin{equation} \label{under}
\xymatrix@R0cm{A[T^{\frac 1S}]_{\underline \sigma} \ar[r] & \mathrm D_{A/R, \underline \sigma}
\\ T^{\frac 1n} \ar@{|->}[r] & 1 + (\sigma_{n}(x) -x) \partial_{\sigma_{n}}}
\end{equation}
inducing a functor
\begin{equation} \label{funcev}
\partial_{A,\underline \sigma}\mathrm{-Mod} \to \underline \sigma_{A}\mathrm{-Mod}.
\end{equation}
If $x$ is a strong $\underline q$-coordinate, then the map \eqref{under} is an isomorphism and the functor \eqref{funcev} is an equivalence.
\end{prop}

\begin{pf}
This will follow from theorem 6.13 of \cite{LeStumQuiros15a} once we know that the various maps for different $n$'s are all compatible.
More precisely, we have to check that the differential operators $1 + (\sigma_{n}(x) -x) \partial_{\sigma_{n}}$ form a system of roots in $\mathrm D_{\underline \sigma}$.
But this follows from proposition \ref{compuiss}.
$\quad \Box$
\end{pf}

\begin{thm}[Formal quantum confluence 2] \label{mainth}
Let $R$ be a $\mathbb Q$-algebra, $\underline q$ a strongly admissible system of roots in $R$ and $A$ a $\underline q$-$R$-algebra.
Then, there exists a canonical $A$-linear map  
\begin{equation*}  \label{fqc2} 
\mathrm D_{A/R, \underline \sigma} \to \widehat {\mathrm D}_{A/R}.
\end{equation*}
If $R$ is $\underline q$-divisible and $\underline q$ is infinite, then the map \eqref{confilm} has dense image.
\end{thm}

\begin{pf}
First of all, the map \eqref{confilm} is obtained by taking the direct limit of the maps \eqref{confone} for all $\sigma_{n}$'s.
Now, we set $p_{n} := q_{n}\mathrm{-char}(R)$ for all $n \in S$ and we apply the confluence theorem \ref{confalg} to $\sigma_{n}$.
If $p_{n} = 0$, we are done. Otherwise, the theorem tells us that the bottom map in the following commutative diagram is surjective:
\begin{displaymath}
\xymatrix{
\mathrm D_{\underline \sigma} \ar[r] & \widehat {\mathrm D} \ar[d]
\\ \mathrm D_{\sigma_{n}} \ar@{->>}[r] \ar[u] & {\mathrm D}/\partial^{p_{n}}
}.
\end{displaymath}
And we proved in lemma \ref{charcri} that for all $n \in S$, we have $p_{n} > n$.
Since $\underline q$ is infinite, we have $p_{n} \to \infty$ and we see that the image of the upper map is dense.
$\quad \Box$
\end{pf}

\begin{cor}
Assume moreover, that $x$ is a strong $\underline q$-coordinate on $A$.
Then, the $A$-linear map
\begin{equation*}
\xymatrix@R0cm{A[T^{\frac 1S}]_{\underline \sigma} \ar[r] & \widehat{\mathrm D}_{A/R}
\\ T^{\frac mn} \ar@{|->}[r] & \sum_{k=0}^{\infty} \frac 1{k!}(\sigma_{n}^m(x) -x)^{k} \partial^{k}}
\end{equation*}
has dense image.
\end{cor}

\begin{pf}
The formula comes from corollary \ref{corsig}
$\quad \Box$
\end{pf}

\begin{xmp}
Theorem \ref{mainth} (and its corollary) applies in particular when $q$ is a primitive $p$-th root of unity in some algebraically closed field $K$ of characteristic zero, $\underline q$ is a system of $p^n$-th roots of $q$ in $K$, $R$ is an algebra containing $K$, $A = R[x, x^{-1}]$ and $\sigma_{n}(x) = q_{n}x$ (we should actually write $\sigma_{p^n}$ and $q_{p^n}$ but we will try to make the notations easier to use).

As we already saw, in this situation, $\mathrm D := \mathrm D_{A/R}$ is the non commutative ring $R[x, x^{-1}, \partial]$ with the commutation rule $\partial x = x \partial +1$ and we may see $\widehat{\mathrm D}$ (which is \emph{not} a ring) as $R[x, x^{-1}][[\partial]]$.
We also want to understand the left hand side and see that $\mathrm D_{n} := \mathrm D_{A/R,\sigma_{n}}$ is the non commutative ring $R[x, x^{-1}, \partial_{n}]$ with the commutation rule $\partial_{n} x = q_{n}x \partial_{n} +1$.
The transition maps are given by the rather tricky formulas
\begin{displaymath}
\partial_{n} \mapsto \sum_{k=1}^{p^{n'-n}} \frac { q_{n'}^{\frac {k(k-1)}{2}}(q_{n'} - 1)^{k-1}} {(k)_{q_{n'}}} {p^{n'-n}-1 \choose k-1}_{q_{n'}}x^{k-1}  \partial_{n'}^{k},
\end{displaymath}
and the map of the theorem satisfies
\begin{displaymath}
\partial_{n} \mapsto \sum_{k=1}^{\infty} \frac 1 {k!}(q_{n}-1)^{k-1}x^{k-1} \partial^{k}.
\end{displaymath}
Also, the map of the corollary is given by
\begin{displaymath}
T^{r} \mapsto \sum_{k=1}^{\infty} \frac 1 {k!}(q^{r} -1)^{k} x^{k}\partial^{k}
\end{displaymath}
if we set $q^{m/n} := q_{n}^m$.
\end{xmp}

\begin{rmk}
Alexei Belov-Kanel and Maxime Konsevich proved in \cite{BelovKanelKonsevich07} that the Jacobian conjecture is stably equivalent to the Dixmier conjecture  \cite{Dixmier68}.
According to the later, any endomorphism of a Weyl algebra over $\mathbb C$ is an automorphism.
And this is still a conjecture even in dimension one.
On the other hand, when $q$ is a primitive $p$-th root of unity, then the quantum Weyl algebra $D_{\sigma}$ does \emph{not} satisfy the Dixmier conjecture but it is an Azumaya algebra.
In particular, checking that an endomorphism is an automorphism can be done on the center.
As explained by Backelin in \cite{Backelin11}, it is appealing to attack the Dixmier conjecture through quantum Weyl algebras and one can hope that theorem \ref{mainth} might provide a tool for this quest. 
\end{rmk}


\bibliographystyle{plain}
\addcontentsline{toc}{section}{References}
\bibliography{BiblioBLS}

\end{document}